\documentclass[a4paper,10pt]{article}
\title{PBW-deformations and deformations \emph{\`a la Gerstenhaber} of $N$-Koszul algebras}
\author{Estanislao Herscovich, Andrea Solotar and Mariano Su\'arez-\'Alvarez
\thanks{\footnotesize This work has been supported by the projects  UBACYTX212 and 475, PIP-CONICET 112-200801-00487, PICT-2007-02182 and MATHAMSUD-NOCOSETA.
The first author is a post-doctoral fellow of the Alexander von Humboldt Foundation.
The second and third authors are  research members of CONICET (Argentina). 
}}
\date{}

\usepackage{amsmath,amsthm,amsfonts,amssymb,stmaryrd}
\usepackage{latexsym}
\usepackage{color}
\usepackage{graphicx}
\usepackage{float}  
\usepackage[small, bf, margin=90pt, tableposition=bottom]{caption}
\usepackage[T1]{fontenc}
\usepackage{palatino}
\usepackage[fleqn,tbtags]{mathtools}

\usepackage[
  breaklinks,
  naturalnames,
  ]{hyperref}

\usepackage{amsrefs}

\input{xy}
\xyoption{all}
\xyoption{poly}
\usepackage[all]{xy}

\newtheorem{theorem}{Theorem}[section]
\newtheorem{theorem*}{Theorem}

\newtheorem{lemma}[theorem]{Lemma}
\newtheorem{proposition}[theorem]{Proposition}
\newtheorem{definition}[theorem]{Definition}
\newtheorem{remark}[theorem]{Remark}
\newtheorem{example}[theorem]{Example}

\newtheorem{fact}[theorem]{Fact}

\numberwithin{equation}{section}                    

\let\oldqed\qed
\renewcommand\qed{\oldqed\par\bigskip}

\newcommand\cl[1]{{\langle#1\rangle}}

\newcommand\CC{{\mathbb{C}}}

\newcommand\ZZ{{\mathbb{Z}}}
\newcommand\NN{{\mathbb{N}}}

\def\Z{{\mathcal Z}}

\def\Tor{{\mathrm {Tor}}}

\def\place{{-}}

\begin{document}

\maketitle

\begin{abstract}
In this article we establish an explicit link between the classical theory of deformations \textit{\`a la Gerstenhaber} 
(and \textit{a fortiori} with the Hochschild cohomology) and (weak) PBW-deformations of homogeneous algebras. 
Our point of view is of cohomological nature.
As a consequence, we recover a theorem by R. Berger and V. Ginzburg, which gives a precise condition for a filtered algebra 
to satisfy the so-called \emph{PBW property}, under certain assumptions. 
\end{abstract}

\textbf{2000 Mathematics Subject Classification:} 16E40, 16S37, 16S80. 

\textbf{Keywords:} Deformation theory, Koszul algebras, Hochschild cohomology.

\section*{Introduction}

Given a graded $k$-algebra $A= TV/\cl{R}$ with $R \subseteq V^{\otimes N}$ there are two available notions of deformation 
of $A$: PBW-deformations and weak PBW-deformations, as defined by R. Berger and V. Ginzburg in \cite{BG06} 
and classical deformations (after M. Gerstenhaber). 
The main goal of this article is to construct explicit equivalences between these concepts, under suitable hypotheses.
Our construction is strongly related to Hochschild cohomology theory.
One of our main motivations is to study the deformation theory of several examples of (graded) $N$-Koszul algebras of interest, 
and since the Hochschild cohomology of many of these algebras is known, 
we believe that it is quite fruitful to have such a direct connection. 

We shall briefly explain our results in more detail. 
We consider a semisimple ring $k$ containing a field $F$ of characteristic zero, such that $k^{e} = k \otimes_{F} k^{op}$ is semisimple. 
Let $A = TV/\cl{R}$ be an $N$-homogeneous $k$-algebra (\textit{i.e.} $R \subseteq V^{\otimes N}$). 
We are interested in studying filtered algebras $U = TV/\cl{P}$ with $P \subseteq \oplus_{i=0}^{N} V^{\otimes i}$, 
such that $R = \pi_{N}(P)$, for $\pi_{N} : TV \rightarrow V^{\otimes N}$ the canonical projection, which satisfy an extra property: 
the surjective morphism of graded algebras $p : A \rightarrow \mathrm{gr}(U)$ induced by the projection $TV \rightarrow U$ is an isomorphism. 
In this case $U$ is called a \textit{PBW-deformation} of $A$. 
The condition of $p$ being injective is equivalent to an infinite number of equalities between certain $k$-bimodules given by intersections 
of tensor powers of $V$ and $R$ (see Subsection \ref{subsec:pbw}). 
The filtered algebra $U$ is said to be a \textit{weak PBW-deformation} of $A$ if only two of this set of equalities are satisfied (see \eqref{eq:pbw1} and \eqref{eq:pbw2}). 

We provide explicit constructions under the hypothesis that $\Tor^{A}_{3}(k,k)$ is concentrated in degree $N+1$, as follows: 
given a graded deformation of $A$ we construct a weak PBW-deformation of $A$, and, conversely, given a weak PBW-deformation of $A$ 
we construct a graded deformation, such that both constructions are inverse up to equivalence. 
They are completely explicit (see Subsections \ref{subsec:pbwtoger} and \ref{subsec:pbwtoger}). 
On the other hand, it is well-known, without the assumption on the torsion group of $A$, that for any filtered algebra $U$ satisfying the PBW property 
one can naturally obtain a graded deformation $A_{t}$ of $A$ by considering the Rees algebra $R(U)$, 
and conversely, given a graded deformation $A_{t}$ of $A$, one obtains a filtered algebra $U$ satisfying the PBW property taking a generic fiber $A_{t}/\cl{t-1}$. 
These two constructions provide quasi-inverse functors for the respective categories. 
Finally, under the assumption that $\Tor^{A}_{3}(k,k)$ is concentrated in degree $N+1$, we notice that the explicit constructions we defined 
in fact coincide with the previous general ones, implying Theorem 1.2 and Theorem 3.4 of \cite{BG06}, which state that 
the notion of a weak PBW-deformation and of a PBW-deformation coincide 
if $\Tor^{A}_{3}(k,k)$ is concentrated in degree $N+1$. 
We would like to point out that the procedure exhibited here allows to find the deformed product of $A_{t}$ explicitly from the filtered algebra, 
even though the computations are often very hard to perform in general. 

We would also like to remark that a similar construction has been already considered in the (second part of the) proof of Theorem 1.1 of the article \cite{FV06} of G. Fl{\o}ystad and J.E. Vatne. 
However, the mentioned proof contains a mistake. 
More precisely, following the notation of that article, on page 122, after the identity defining $\gamma$ on line 22, it is stated that $\gamma \circ \sigma$ 
vanishes. 
This is not necessarily true, because there is in principle no identification of the Koszul resolution $K_{\bullet}$ inside the bar resolution $B_{\bullet}$, which is compatible with taking brackets $[1, \place]$. 
In fact, this can also be noticed from the fact that the recursion formulas for the cochains giving the deformed product corresponding to a weak PBW-deformation explained in Proposition \ref{prop:gertopbw} of 
Subsection \ref{subsec:pbwtoger} are indeed more complicated. 

The article is organized as follows. 
In the first section we recall some generalities about the theory of $N$-homogeneous algebras over a (not necessarily commutative) semisimple ring 
$k$ containing a field $F$ of characteristic zero, such that $k^{e} = k \otimes_{F} k^{op}$ is semisimple, 
and make an intensive study of the (reduced) Hochschild resolution and a ``Koszul-like'' projective resolution of bimodules of an 
$N$-homogeneous algebra satisfying some vanishing condition on the torsion groups. 
In Section \ref{sec:pbwger}, we first recall general facts on PBW-deformations and weak PBW-deformations of $N$-homogeneous algebras 
and the classical graded deformation theory \textit{\`a la Gerstenhaber}. 
Later, in Subsections \ref{subsec:gertopbw} and \ref{subsec:pbwtoger} we establish a link between these concepts. 
Finally, in Section \ref{sec:main} we prove the main results of this paper, namely Theorems \ref{teo:pbw}, \ref{teo:pbw+} and \ref{teo:final} and we give several examples. 

\section*{Acknowledgements}

We would like to thank Roland Berger for several comments which helped to improve the article. 
We would also like to thank the referee for a careful reading of the manuscript. 

\section{Generalities}
\label{sec:gen}

From now on, $k$ will be a (not necessarily commutative) semisimple ring containing a field $F$ of characteristic zero. 
We consider $k^{e} = k \otimes_{F} k^{op}$. 
We assume further that $k^{e}$ is semisimple. 
By \emph{$k$-bimodule} we will always mean a $k$-bimodule such that the action of $F$ is symmetric, \textit{i.e.} $a.m = m.a$, for all $a \in F$, and $m$ in the $k$-bimodule $M$. 
Note that we do not assume the action of $k$ on the bimodule to be symmetric. 
It is clear that this definition of $k$-bimodule is equivalent to the notion of a (say left) $k^{e}$-module. 
As a matter of notation, all unadorned tensor products $\otimes$ are over $k$.

A \emph{$k$-algebra} denotes a monoid object in the monoidal category of $k$-bimodules, \textit{i.e.} 
it is a $k$-bimodule $A$ provided with a morphism $\mu : A \otimes A \rightarrow A$ of $k$-bimodules, which will be denoted $\mu(a \otimes a') = a.a'$, 
and an element $1 \in A$, such that $\mu$ is associative, $1.a = a.1 = a$, for all $a \in A$, and $c.1 = 1.c$, for all $c \in k$. 
Equivalently, a $k$-algebra is a ring $A$ provided with a unitary ring homomorphism $i_{A} : k \rightarrow A$ where $\mathrm{Im}(i_{A})$ is not necessarily contained in $\Z(A)$. 
For $n \geq 2$, $\mu^{(n)} : A^{\otimes n} \rightarrow A$ will denote the morphism of $k$-bimodules defined recursively by 
$\mu^{(2)} = \mu$ and $\mu^{(n+1)} = \mu^{(n)} \circ (\mu \otimes 1_{A}^{\otimes (n-1)})$, for $n > 2$.  
A \emph{morphism} from a $k$-algebra $A$ to a $k$-algebra $B$ is 
a ring homomorphism $f : A \rightarrow B$ such that $f \circ i_{A} = i_{B}$.

A \emph{graded $k$-bimodule} will be a $k$-bimodule $V$ together with a decomposition as a direct sum of $k$-bimodules 
$V = \oplus_{n \in \ZZ} V_{n}$.  
A \emph{morphism} of graded $k$-bimodules is just a degree preserving morphism of $k$-bimodules. 
The category of graded $k$-bimodules is monoidal in the obvious manner. 
We will denote the homomorphism group between two graded $k$-bimodules $M$ and $M'$ by $\mathrm{hom}_{k^{e}}(M,M')$. 
The \emph{shift functor} $(\place)[1]$, together with its iterations, is defined in the usual way, and we recall that the internal group of homomorphisms between two graded $k$-bimodules 
$M$ and $M'$ is given by $\mathcal{H}om_{k^{e}}(M,M') = \oplus_{i \in \ZZ} \mathrm{hom}_{k^{e}}(M,M'[i])$.
Moreover, a \emph{graded $k$-algebra} is a monoid object in the monoidal category of graded $k$-bimodules, 
\textit{i.e.} a $k$-algebra provided with a decomposition of $k$-bimodules of the form $A = \oplus_{n \in \ZZ} A_{n}$ such that 
$1 \in A_{0}$ and $A_{n}.A_{m} \subseteq A_{n+m}$, for all $n, m \in \ZZ$. 
In fact, we shall usually assume that $A = \oplus_{n \in \NN_{0}} A_{n}$. 

Given a graded $k$-algebra $A$, a \emph{graded left $A$-module} $M$ is given by a graded $k$-bimodule structure on $M$ together with a morphism of graded $k$-bimodules 
$\rho : A \otimes M \rightarrow M$, which will be denoted by $\rho(a \otimes m \otimes a') = a.m$, for $a \in A$ and $m \in M$ satisfying the usual mixed associative axiom, \textit{i.e.} 
$a.(a'.m) = (a.a').m$, for all $a, a' \in A$ and $m \in M$, and that $1.m = m$, for all $m \in M$. 
As usual, a \emph{morphism} of graded left $A$-modules is just a degree preserving morphism of $A$-modules, and 
the homomorphism group between two graded left $A$-modules $M$ and $M'$ will be denoted by $\mathrm{hom}_{A}(M,M')$. 
We note that shift functor $(\place)[1]$ may also be defined on a graded left $A$-module $M$, where the underlying structure of graded $k$-bimodule is the same as before, 
and the obvious action. 
As in the previous case, the internal group of homomorphisms between two graded left $A$-modules $M$ and $M'$ is given by 
$\mathcal{H}om_{A}(M,M') = \oplus_{i \in \ZZ} \mathrm{hom}_{A}(M,M'[i])$.
The usual definition of \emph{graded right $A$-module} and \emph{graded $A$-bimodule} are analogous, 
and using the obvious structure of graded $k$-algebra on $A^{e} = A \otimes A^{op}$ we also see that the notion graded $A$-bimodule and graded (say left) $A^{e}$-module coincide. 

Let $N$ be a positive integer, $N \geq 2$. 
By \emph{complex} (resp. $N$-\emph{complex}) we mean a nonnegatively graded module over a $k$-algebra $A$, provided with an endomorphism $d$ of degree $-1$ 
such that $d^{2}=0$ (resp. $d^{N}=0$). 

We note that, since $k$ is semisimple, it is von Neumann regular, so all the considerations in~\cite{BG06}, Section 2, in order to properly consider 
the notion of Koszul algebra also apply to this case. 
Let $A = TV/\cl{R}$ be an \emph{$N$-homogeneous algebra}, where $V$ is a 
$k$-bimodule (considered to be concentrated in degree $1$) 
and $R$ is a $k$-subbimodule of $V^{\otimes N}$. 
In this situation, we shall identify $A/k$ with the $k$-subbimodule $I_{+}$ of $A$ spanned by elements of strictly positive degree. 
We will make use of the number function given by
\begin{align*}
  \zeta : \NN_{0} &\rightarrow \NN_{0}
   \\
   \zeta(n) &= \begin{cases}             
                     N m, &\text{if $n = 2 m$,}
                     \\
                     N m+1, &\text{if $n = 2 m +1$.}
                    \end{cases}
\end{align*}
For $n \in \NN_{0}$, define $W_{n}$ the $k$-subbimodule of $V^{\otimes n}$ given by $V^{\otimes n}$ if $n < N$, and by 
\[     \bigcap_{i=0}^{n-N} V^{\otimes i} \otimes R \otimes V^{\otimes (n-i-N)},     \]
if $n \geq N$.
Then $W_{i} = V^{\otimes i}$, for $i=0,\dots,N-1$, $W_{N} = R$, and 
$W_{N+1} = (R \otimes V) \cap (V \otimes R)$, which we shall also denote by $R_{N+1}$. 
Note that $W_{n}$ may be regarded as a graded $k$-bimodule concentrated in degree $n$. 

We recall the \emph{bimodule Koszul complex} $(K_{\bullet}(A),d_{\bullet})$ of $A$, defined in~\cite{BM06} for the case that $k = F$ is a field. 
First, we consider the graded $A$-bimodule given by $(A \otimes W_{n} \otimes A)_{n \in \NN_{0}}$. 
Then, for each $n \in \NN$, there are two $A$-bimodule maps
\[     d^{L}_{n}, d^{R}_{n} : A \otimes W_{n} \otimes A \rightarrow A \otimes W_{n-1} \otimes A     \]
given by 
\begin{align*}
   d^{L}_{n}(a_{0} \otimes \dots \otimes a_{n+1}) &= a_{0} a_{1} \otimes a_{2} \otimes \dots \otimes a_{n+1}, 
   \\
   d^{R}_{n}(a_{0} \otimes \dots \otimes a_{n+1}) &= a_{0} \otimes \dots \otimes a_{n-1}\otimes a_{n} a_{n+1}.     
\end{align*}
We note that $d^{L}$ and $d^{R}$ commute. 
We shall denote an element $a \otimes \alpha \otimes b \in A \otimes W_{n} \otimes A$, where $a,b \in A$ and $\alpha \in W_{n}$, 
in the shorter form $a|\alpha|b$. 

If $q$ is a primitive $N$-th root of unity, we consider the map of $A$-bimodules 
\[     d^{b}_{n} : A \otimes W_{n} \otimes A \rightarrow A \otimes W_{n-1} \otimes A     \]
given by $d^{b}_{n} = d^{L}_{n} - q^{n-1} d^{R}_{n}$. 
It is trivially verified that $((A \otimes W_{\bullet} \otimes A),d^{b}_{\bullet})$ is an $N$-complex. 
The bimodule Koszul complex $(K_{\bullet}(A),d_{\bullet})$ of $A$ is the $(1,N-1)$-contraction of the previous $N$-complex, that is $K_{n}(A) = A \otimes W_{\zeta(n)} \otimes A$ 
with differential $d_{\bullet}$ given by the corresponding successive composition of the differential of the previous $N$-complex. 
It is easy to see that 
\begin{equation}
\label{eq:differential}
     d_{n} = \begin{cases}
                      d^{L}_{N m +1} - d^{R}_{N m + 1}, &\text{if $n = 2 m + 1$,}
                      \\
                      d^{L}_{N (m - 1) + 2} \dots d^{L}_{N m} 
                     + d^{R}_{N (m - 1) + 2} \dots d^{R}_{N m}  &
                     \\+ \sum\limits_{i=0}^{N-3} d^{L}_{N (m - 1) + 2} \dots d^{L}_{N m - i -1} d^{R}_{N m - i} \dots d^{R}_{N m}, &\text{if $n = 2 m$.}
                      \end{cases}                   
\end{equation}
We thus notice that the $N$-th root of unity is in some sense superfluous, since we may define the Koszul complex without invoking it. 
The algebra $A$ is called \emph{Koszul} if the bimodule Koszul complex $(K_{\bullet}(A),d_{\bullet})$ is a resolution of $A$-bimodules of $A$ for the map 
$d_{0} = \mu : K_{0}(A) = A \otimes A \rightarrow A$ given by the product 
$\mu$ of $A$. 
We recall that the category of $A$-bimodules is equivalent to the category of 
(say left) modules over 
$A^{e} = A \otimes_{F} A^{op}$ (not over $A^{e} = A \otimes_{k} A^{op}$, because we are not considering symmetric $k$-bimodules!). 
Moreover, since $\mathcal{H}om_{A^{e}} (A \otimes W_{n} \otimes A , \place) \simeq \mathcal{H}om_{k^{e}}(W_{n},\place)$, 
we see that the bimodule Koszul resolution consists of projective graded $A$-bimodules, \textit{i.e.} projective graded left $A^{e}$-modules. 

For our purposes, we will be interested in a weaker condition than that of being Koszul: 
we shall suppose that $A$ is an $N$-homogeneous algebra satisfying that $\Tor_{3}^{A}(k,k)$ is concentrated in degree $N+1$. 
This is equivalent to say that there exists a projective resolution (of graded $A$-bimodules) of $A$ that coincides with the one given 
by the Koszul bimodule complex for homological degrees less than or equal to three (see \cite{BG06}, Sec. 2). 
This assumption comes from the fact that all the computations we shall perform in this article are restricted to those homological degrees. 
In this case, we shall still denote by $K_{\bullet}(A)$ the former projective resolution (of graded $A$-bimodules) of $A$, and call it \emph{minimal}.  

The \emph{(graded) Hochschild complex} $(C_{\bullet}(A),b_{\bullet})$ of $A$ is given by 
$C_{n}(A) = A \otimes A^{\otimes n} \otimes A$, for $n \in \NN_{0}$, with differential 
\[     b_{n} (a_{0} \otimes \dots \otimes a_{n+1}) = \sum_{i=0}^{n} (-1)^{i} a_{0} \otimes \dots \otimes a_{i} a_{i+1} \otimes \dots \otimes a_{n+1},     \]
for $n \in \NN$.
We will also write $a_{0}|\dots|a_{n+1}$ instead of $a_{0} \otimes \dots \otimes a_{n+1}$. 
This complex gives a projective resolution $C_{\bullet}(A) \rightarrow A$ of $A$ in the category of (graded) $A$-bimodules provided with morphisms of degree zero 
via $b_{0} = \mu : A \otimes A \rightarrow A$, called the \emph{Hochschild resolution}. 
Moreover, we may consider the \emph{reduced (graded) Hochschild resolution} $(\bar{C}_{\bullet}(A), \bar{b}_{\bullet})$ of $A$. 
The underlying graded $A$-bimodule is given by $\bar{C}_{n}(A) = A \otimes (A/k)^{\otimes n} \otimes A$, for $n \in \NN_{0}$, so there 
exists a canonical projection $p_{n} : C_{n}(A) \rightarrow \bar{C}_{n}(A)$, and the differential $\bar{b}_{n}$ is induced by $b_{n}$, for $n \in \NN$, such that 
$p_{\bullet}$ becomes a morphism of complexes. 
It gives a projective resolution of $A$ in the category of (graded) $A$-bimodules provided with morphisms of degree zero, 
called the \emph{reduced Hochschild resolution} $\bar{C}_{\bullet}(A) \rightarrow A$ of $A$. 
Using the identification between $A/k$ and the $k$-subbimodule of $A$ spanned by the strictly positive elements of $A$, 
we will usually denote an element of $\bar{C}_{n}(A)$ simply by $a_{0} \otimes a_{1} \otimes \dots \otimes a_{n} \otimes a_{n+1}$, 
or $a_{0}|a_{1}|\dots|a_{n}|a_{n+1}$, where $\deg(a_{i}) \geq 1$, for $i=1,\dots,n$, 
instead of the more correct $a_{0} \otimes [a_{1}] \otimes \dots \otimes [a_{n}] \otimes a_{n+1}$ 
(or $a_{0}|[a_{1}]|\dots|[a_{n}]|a_{n+1}$), where $[a] \in A/k$ is the class of $a \in A$. 
Furthermore, if by chance there exist $i \in \{ 1, \dots, n \}$ such that $\deg(a_{i}) = 0$, we may also consider $a_{0}|a_{1}|\dots|a_{n}|a_{n+1}$ as the zero element of 
$\bar{C}_{n}(A)$. 
We remark that, under this identification, $\bar{b}_{n}$ is the restriction of $b_{n}$ to $A \otimes I_{+}^{\otimes n} \otimes A$. 


If $M$ is a graded $A$-bimodule, the \emph{graded Hochschild cohomology groups} of $A$ with coefficients in $M$, which will be also denoted by $H^{\bullet}(A,M)$, are given by the cohomology of the cochain complex 
\begin{multline*}
     \mathcal{H}om_{A^{e}}(C_{\bullet}(A),M) = \bigoplus_{j \in \ZZ} \mathrm{hom}_{A^{e}}(C_{\bullet}(A),M[j])  
     \\
     \simeq \mathcal{H}om_{k^{e}}(A^{\otimes \bullet},M) = \bigoplus_{j \in \ZZ} \mathrm{hom}_{k^{e}}(A^{\otimes \bullet},M[j])     
\end{multline*}
with the induced differential, and where $\mathrm{hom}(\place,\place)$ is the set of degree preserving homomorphisms. 
Again, these cohomology groups can be computed using either the complex $\mathcal{H}om_{A^{e}}(C_{\bullet}(A),M)$ or 
$\mathcal{H}om_{A^{e}}(\bar{C}_{\bullet}(A),M)$. 
If $M = A$ again, we shall write $HH^{\bullet}(A)$ instead of $H^{\bullet}(A,A)$. 
In this case, we know that $H^{\bullet}(A,M) \simeq \mathcal{E}xt_{A^{e}/k^{e}}^{\bullet} (A,M)$, the relative derived functors of $\mathcal{H}om_{A^{e}}(A,\place)$. 
We see that $H^{\bullet}(A,M)$ has an internal grading or weight that comes from the gradings of $A$ and $M$, which we shall denote by $H^{\bullet}(A,M)_{\bullet}$. 

Before proceeding further, we shall state some notation. 
Let us consider a positive integer $p \leq n$ and a subset $I = \{ i_{1} < \dots < i_{m} \}$ of $\{ 1, \dots, p \}$. 
Given a decomposition $n = n_{1} + \dots + n_{p}$ of $n$ ($n_{j} \in \NN_{0}$) satisfying that $n_{i_{j}} = 1$, for $j = 1, \dots, m$, 
we define the collection of $k^{e}$-linear maps $q_{i}^{(n_{1},\dots,n_{p})}$ homogeneous of degree zero given by the canonical map
$q_{i}^{(n_{1},\dots,n_{p})} : V^{\otimes n_{i}} \rightarrow A$, if $i \notin I$, and $q^{(n_{1},\dots,n_{p})}_{i} = 1_{V}$, if $i \in I$. 
Then, consider the map $\mathrm{par}_{(n_{1},\dots,n_{p}),I} = q^{(n_{1},\dots,n_{p})}_{1} \otimes \dots \otimes q^{(n_{1},\dots,n_{p})}_{p}$, 
so it is a map from $V^{\otimes n}$ to 
\[     A^{\otimes (i_{1}-1)} \otimes V \otimes A^{\otimes (i_{2} - i_{1} - 1)} \otimes V \otimes 
       \dots \otimes V 
       \otimes A^{\otimes (i_{m} - i_{m-1} - 1)} \otimes V \otimes A^{\otimes (p - i_{m})}.     \]
We now define the maps $\mathrm{par}'_{(n_{1},\dots,n_{p}),I}$ and $\mathrm{Par}_{p,I}$ from $V^{\otimes n}$ to 
$A \otimes V \otimes A \otimes V \otimes \dots \otimes V \otimes A \otimes V \otimes A$, 
where the $k$-bimodule $V$ appears $m$ times in the last tensor product. 
The former is given by the composition of 
\[     (\mu^{(i_{1}-1)} \otimes 1_{V} \otimes \mu^{\otimes (i_{2} - i_{1} - 1)} \otimes 1_{V} \otimes \dots \otimes 1_{V} 
       \otimes \mu^{(i_{m} - i_{m-1} - 1)} \otimes 1_{V} \otimes \mu^{(p - i_{m})})     \] 
and $\mathrm{par}_{(n_{1},\dots,n_{p}),I}$. 
The latter is defined as $\mathrm{Par}_{p,I} = \sum \mathrm{par}'_{(n_{1},\dots,n_{p}),I}$, 
where the last sum is indexed over all decompositions $(n_{1},\dots,n_{p})$ of $n$ such that $n_{i_{j}} = 1$, for $j = 1, \dots, m$. 
In general, we shall write 
\[     \mathrm{Par}_{p,I}(\alpha) 
       = \alpha_{(1)} \otimes \bar{\alpha}_{(2)} \otimes \alpha_{(3)} \otimes \dots \otimes 
       \alpha_{(2m-1)} \otimes \bar{\alpha}_{(2 m)} \otimes \alpha_{(2m+1)},     \]
where the bars are used for the elements with indices in $I$. 
We remark that in the previous expression a sum over all decompositions of $\alpha$ in $p$ terms such that the ones with bars belong to $V$ is implicit. 
For instance, if $\alpha = v_{1} \dots v_{n} \in V^{\otimes n}$, we denote 
\[     \mathrm{Par}_{3,\{2\}}(\alpha)=  \alpha_{(1)}|\bar{\alpha}_{(2)}|\alpha_{(3)} \in A \otimes V \otimes A     \] 
the sum over all decompositions of $\alpha$ such that the term with the bar $\bar{\alpha}_{(2)}$ belongs to $V$, 
\textit{i.e.} it is given by 
\[     1|v_{1}|v_{2} \dots v_{n} + \sum_{i=1}^{n-2} v_{1} \dots v_{i} | v_{i+1} | v_{i+2} \dots v_{n} + v_{1} \dots v_{n-1} | v_{n} | 1.     \]
This should be seen as a similar notation to Sweedler's one for coproducts.
We emphasize that each term $\alpha_{(i)}$ is homogeneous and may be (and shall be) taken equal to $1$ if it has degree $0$. 

The following remarks are easy consequence of the definitions:
\begin{fact}
\label{fac:cuenta division}
If $\alpha = v_{1} \dots v_{n} \in V^{\otimes n}$ and $\beta = w_{1} \dots w_{m} \in V^{\otimes m}$, 
then 
\[     (\alpha \beta)_{(1)}|\overline{(\alpha \beta)}_{(2)}|(\alpha \beta)_{(3)} 
       = \alpha_{(1)}|\bar{\alpha}_{(2)}|\alpha_{(3)} \beta + \alpha \beta_{(1)}|\bar{\beta}_{(2)}|\beta_{(3)}.     \]
\end{fact}

\begin{fact}
\label{fac:cuenta telescopica}
If $\alpha = v_{1} \dots v_{n} \in V^{\otimes n}$, let us consider $\alpha_{(1)}|\bar{\alpha}_{(2)}|\alpha_{(3)} \in A \otimes V \otimes A \subseteq 
 A \otimes A \otimes A$. 
Then, 
\[     b_{1} (\alpha_{(1)}|\bar{\alpha}_{(2)}|\alpha_{(3)}) = 
       \alpha_{(1)} \bar{\alpha}_{(2)}|\alpha_{(3)} - \alpha_{(1)}|\bar{\alpha}_{(2)}\alpha_{(3)} 
       = \alpha | 1 - 1 | \alpha.     \]
\end{fact}

We see that 
\begin{align}  
     d_{2} (1|r|1) &= r_{(1)} | \bar{r}_{(2)} | r_{(3)},   
     \label{eq: dif r}
     \\    
     d_{3} (1|w|1) &= v_{i}| s_{i} | 1 - 1 | r_{i} | u_{i} = \bar{w}_{(1)}| w_{(2)} | 1 - 1 | w_{(1)} | \bar{w}_{(2)},   
     \label{eq: dif w}
\end{align}
where $r \in R$, and $w = \sum_{i \in I} r_{i} u_{i} = \sum_{i \in I} v_{i} s_{i} \in (V \otimes R) \cap (R \otimes V) = R_{N+1}$, for $u_{i}, v_{i} \in V$, $r_{i},s_{i} \in R$ and a finite set $I$ of indices. 
We have omitted the sum in \eqref{eq: dif w} and we shall do so for the typical elements of $R_{N+1}$: 
we shall simply write $w = r_{i} u_{i} = v_{i} s_{i}$. 
Moreover, the previous choice of letters will be the usual one to denote elements of $R$ and $R_{N+1}$, unless we say the contrary. 

Since both the Hochschild and the minimal bimodule resolutions are projective resolutions of graded $A$-bimodules of $A$, 
there exists unique (up to homotopy) comparison morphisms of graded $A$-bimodules $\sigma_{\bullet} : K_{\bullet}(A) \rightarrow C_{\bullet}(A)$ and 
$\tau_{\bullet} : C_{\bullet}(A) \rightarrow K_{\bullet}(A)$. 
We define the morphisms of graded $A$-bimodules given by the extension of the following expressions 
\begin{align}
   &\sigma_{0} = 1_{A} \otimes 1_{A}, 
   \label{eq:sigma0}\addtocounter{equation}{1}\tag{SIGMA$_{0}$}
   \\
   &\sigma_{1} = 1_{A} \otimes \mathrm{inc}_{V \subset A} \otimes 1_{A},
   \label{eq:sigma1}\addtocounter{equation}{1}\tag{SIGMA$_{1}$}
   \\
   &\sigma_{2}(1|r|1) = 1|r_{(1)}|\bar{r}_{(2)}|r_{(3)},
   \label{eq:sigma2}\addtocounter{equation}{1}\tag{SIGMA$_{2}$}
   \\
   &\sigma_{3}(1|w|1) = 1|v_{i}|s_{i,(1)}|\bar{s}_{i,(2)}|s_{i,(3)} = 1|\bar{w}_{(1)}|w_{(2)}|\bar{w}_{(3)}|w_{(4)}.
   \label{eq:sigma3}\addtocounter{equation}{1}\tag{SIGMA$_{3}$}
\end{align}
We define $\bar{\sigma}_{\bullet} = p_{\bullet} \circ \sigma_{\bullet}$. 
It is easy to see that $\sigma_{0} d_{1} = b_{1} \sigma_{1}$ and $\sigma_{1} d_{2} = b_{2} \sigma_{2}$. 
These identities immediately imply that $\bar{\sigma}_{0} d_{1} = \bar{b}_{1} \bar{\sigma}_{1}$ and $\bar{\sigma}_{1} d_{2} = \bar{b}_{2} \bar{\sigma}_{2}$ hold. 
We shall check that $\sigma_{2} d_{3} = b_{3} \sigma_{3}$, which also yields $\bar{\sigma}_{2} d_{3} = \bar{b}_{3} \bar{\sigma}_{3}$. 
On the one hand, we see that 
\begin{align*}     
   b_{3} \sigma_{3} (1|w|1) &= b_{3} (1|v_{i}|s_{i,(1)}|\bar{s}_{i,(2)}|s_{i,(3)})
   \\
                            &= v_{i}|s_{i,(1)}|\bar{s}_{i,(2)}|s_{i,(3)} - 1|v_{i} s_{i,(1)}|\bar{s}_{i,(2)}|s_{i,(3)} 
                              + 1|v_{i}|s_{i,(1)} \bar{s}_{i,(2)}|s_{i,(3)} 
   \\                       &       - 1|v_{i}|s_{i,(1)}|\bar{s}_{i,(2)} s_{i,(3)}
   \\
                             &= v_{i}|s_{i,(1)}|\bar{s}_{i,(2)}|s_{i,(3)} - 1|v_{i} s_{i,(1)}|\bar{s}_{i,(2)}|s_{i,(3)} 
                              + 1|v_{i}|s_{i}|1 - 1|v_{i}|1|s_{i}
   \\ 
                             &= v_{i}|s_{i,(1)}|\bar{s}_{i,(2)}|s_{i,(3)} - 1|v_{i} s_{i,(1)}|\bar{s}_{i,(2)}|s_{i,(3)},                            
\end{align*}
where we have used Fact \ref{fac:cuenta telescopica} and that $1|v_{i}|s_{i}|1 = 1|v_{i}|1|s_{i} = 0$ since each $s_{i}$ vanishes in $A$. 
On the other hand, 
\begin{align*}   
   \sigma_{2} d_{3} (1|w|1) &= \sigma_{2} (v_{i}|s_{i}|1 - 1|r_{i}|u_{i}) 
   \\
                            &= v_{i}|s_{i,(1)}|\bar{s}_{i,(2)}|s_{i,(3)} - 1|r_{i,(1)}|\bar{r}_{i,(2)}|r_{i,(3)}u_{i}.
\end{align*}
The equality $\sigma_{2} d_{3} = b_{3} \sigma_{3}$ then follows from the simple result: 
\begin{fact}
\label{fac:cuenta w}
If $w = r_{i} u_{i} = v_{i} s_{i} \in R_{N+1}$, for $u_{i}, v_{i} \in V$, $r_{i},s_{i} \in R$ (summation understood), then 
\[     1|v_{i}s_{i,(1)}|\bar{s}_{i,(2)}|s_{i,(3)} = 1|r_{i,(1)}|\bar{r}_{i,(2)}|r_{i,(3)}u_{i}.     \]
\end{fact}
\noindent\textbf{Proof.}
We point out that in the first member of the previous identity we are considering the decompositions of $w$ in three terms where the first one has degree greater than or equal to $1$ and the second one has degree $1$, 
whereas in the second member we consider the decompositions of $w$ in three terms where the second one has degree $1$ but 
the third one has degree greater than or equal to $1$. 

We may decompose the sum $1|v_{i} s_{i,(1)}|\bar{s}_{i,(2)}|s_{i,(3)}$ in two separate cases: 
when $s_{i,(3)} \in k$ (in which case, it is equal to $1$) and when $s_{i,(3)}$ has degree greater than or equal to $1$. 
We may write this as 
\[     1|v_{i} s_{i,(1)}|\bar{s}_{i,(2)}|s_{i,(3)} =               
       \underset{\deg(s_{i,(3)})>0}{\underbrace{1|v_{i} s_{i,(1)}|\bar{s}_{i,(2)}|s_{i,(3)}}}
       + 1|v_{i} s_{i,(1)}|\bar{s}_{i,(2)}|1.     \]
The same reasoning applies to $1|r_{i,(1)}|\bar{r}_{i,(2)}|r_{i,(3)} u_{i}$, to give
\[     1|r_{i,(1)}|\bar{r}_{i,(2)}|r_{i,(3)} u_{i} =               
       \underset{\deg(r_{i,(1)})>0}{\underbrace{1|r_{i,(1)}|\bar{r}_{i,(2)}|r_{i,(3)} u_{i}}}
       + 1|1|\bar{r}_{i,(1)}|r_{i,(2)} u_{i}.     \]

However, since the terms $1|v_{i} s_{i,(1)}|\bar{s}_{i,(2)}|1$ indicate all the decompositions of $w$ in two terms such that the second one has degree $1$, we see that $1|v_{i} s_{i,(1)}|\bar{s}_{i,(2)}|1 = 1|r_{i}|u_{i}|1$, which vanishes, 
since $1|r_{i}|u_{i}|1 \in C_{2}(A)$ and $r_{i} = 0$ in $A$. 
A similar argument implies that $1|1|\bar{r}_{i,(1)}|r_{i,(2)} u_{i} = 1|1|v_{i}|s_{i} \in C_{2}(A)$ vanishes. 

Hence, we only have to prove that 
\[     \underset{\deg(s_{i,(3)})>0}{\underbrace{1|v_{i} s_{i,(1)}|\bar{s}_{i,(2)}|s_{i,(3)}}} 
       = \underset{\deg(r_{i,(1)})>0}{\underbrace{1|r_{i,(1)}|\bar{r}_{i,(2)}|r_{i,(3)} u_{i}}}.     \]
This identity holds, since both members are built from the decompositions of $w$ in three terms, 
where the first and third ones have degree greater than or equal to $1$, and the second one has degree $1$. 
\qed

\begin{definition}
\label{defi:normalized}
Given an homogeneous element $1|a_{1}|\dots|a_{n}|1 \in C_{n}(A)$ (resp. $a_{1} \otimes \dots \otimes a_{n} \in A^{\otimes n}$), we shall say that it is $N$-\emph{normalized} if $\deg(a_{1}) + \dots + \deg(a_{n}) < N$. 
Since $N$ is fixed throughout this work, in both cases we will more simply say that it is \emph{normalized}. 

If $\sum_{i \in I} 1|a_{1}^{i}|\dots|a_{n}^{i}|1 \in C_{n}(A)$ (resp. $\sum_{i \in I} a_{1}^{i} \otimes \dots \otimes a_{n}^{i} \in A^{\otimes n}$) satisfies that $\deg(a_{1}^{i}) + \dots + \deg(a_{n}^{i}) = N$, $\deg(a_{1}^{i}), \deg(a_{n}^{i}) > 0$ and 
$\sum_{i\in I} a_{1}^{i} \dots a_{n}^{i} = 0$ in $A$ (\textit{i.e.} it belongs to $R$ as an element of the tensor algebra $TV$), we shall say that $\sum_{i \in I} 1|a_{1}^{i}|\dots|a_{n}^{i}|1$ 
(resp. $\sum_{i \in I} a_{1}^{i} \otimes \dots \otimes a_{n}^{i} \in A^{\otimes n}$) is a \emph{relation decomposition}. 
This could be also abbreviated by \emph{rel. decomp}.

Finally, if $\sum_{i \in I} 1|a_{1}^{i}|a_{2}^{i}|a_{3}^{i}|1 \in C_{3}(A)$ 
(resp. $\sum_{i \in I} a_{1}^{i} \otimes a_{2}^{i} \otimes a_{3}^{i} \in A^{\otimes 3}$) satisfies that $\deg(a_{1}^{i}) + \deg(a_{2}^{i}) + \deg(a_{3}^{i}) = N+1$, 
$\deg(a_{1}^{i}), \deg(a_{3}^{i}) = 1$ and $\sum_{i\in I} a_{1}^{i} a_{2}^{i} a_{3}^{i} \in R_{N+1}$ as an element of the tensor algebra $TV$, 
we shall say that $\sum_{i \in I} 1|a_{1}^{i}|a_{2}^{i}|a_{3}^{i}|1$ (resp. $\sum_{i \in I} a_{1}^{i} \otimes a_{2}^{i} \otimes a_{3}^{i} \in A^{\otimes 3}$) 
is a \emph{double relation decomposition}. 
It will be occasionally abbreviated by \emph{double rel. decomp}.

Using the identification $A/k \simeq I_{+}$, we have the analogous versions of the previous three definitions for elements of $\bar{C}_{n}(A)$ or $(A/k)^{\otimes n}$, which will be called in the same way. 
We would like to note that the set of homogeneous elements in $C_{n}(A)$ (resp. $A^{\otimes n}$, $\bar{C}_{n}(A)$, $(A/k)^{\otimes n}$) which are normalized 
form a $k$-subbimodule. 
The same applies to the set of homogeneous elements which are relation decompositions, or double relation decompositions. 
\end{definition}

We partially define the morphisms of $A$-bimodules $\tau_{\bullet}$ as follows:
\begin{small}
\begin{align}
   &\tau_{0} = 1_{A} \otimes 1_{A}, &
   \label{eq:tau0}\addtocounter{equation}{1}\tag{TAU$_{0}$}
   \\
   &\tau_{1}(1|a|1) = a_{(1)}|\bar{a}_{(2)}|a_{(3)}, &\text{if $\deg(a) < N$,}
   \label{eq:tau1}\addtocounter{equation}{1}\tag{TAU$_{1}$}
   \\
   &\tau_{2}(1|a|b|1) = 0, &\text{if $1|a|b|1$ is normalized,}
   \label{eq:tau2,1}\addtocounter{equation}{1}\tag{TAU$_{2,1}$}
   \\
   &\tau_{2}(1|a_{i}|b_{i}|1) = 1|r|1, &\text{if $r = \sum_{i \in I} a_{i} \otimes b_{i}$ is a rel. decomp.,}
   \label{eq:tau2,2}\addtocounter{equation}{1}\tag{TAU$_{2,2}$}
   \\
   &\tau_{3}(1|a|b|c|1) = 0, &\text{if $1|a|b|c|1$ is normalized,}
   \label{eq:tau3,1}\addtocounter{equation}{1}\tag{TAU$_{3,1}$}
   \\
   &\tau_{3}(1|a_{i}|b_{i}|c_{i}|1) = 0, &\text{if $r = \sum_{i \in I} a_{i} \otimes b_{i} \otimes c_{i}$ is a rel. decomp.,}
   \label{eq:tau3,2}\addtocounter{equation}{1}\tag{TAU$_{3,2}$}
   \\
   &\tau_{3}(1|a_{i}|b_{i}|c_{i}|1) = 1|w|1, &\text{if $w = \sum_{i \in I} a_{i} \otimes b_{i} \otimes c_{i}$ is a double rel. decomp.}
   \label{eq:tau3,3}\addtocounter{equation}{1}\tag{TAU$_{3,3}$}
\end{align}
\end{small}
We remark that the previous identities induce maps $\bar{\tau}_{\bullet}$ from $\bar{C}_{\bullet}(A)$ to $K_{\bullet}(A)$ (partially defined on the image under 
$p_{\bullet}$ of the domain of definition of $\tau_{\bullet}$). 
It is trivial to verify that $\tau_{0} b_{1} = d_{1} \tau_{1}$ holds, wherever $\tau_{1}$ is defined. 
We shall check that $\tau_{1} b_{2} = d_{2} \tau_{2}$ and $\tau_{2} b_{3} = d_{3} \tau_{3}$ are verified where we have defined them. 
These identities would imply that $\bar{\tau}_{0} \bar{b}_{1} = d_{1} \bar{\tau}_{1}$, $\bar{\tau}_{1} \bar{b}_{2} = d_{2} \bar{\tau}_{2}$ and 
$\bar{\tau}_{2} \bar{b}_{3} = d_{3} \bar{\tau}_{3}$ hold where they are defined. 
These maps can be extended to complete morphisms of complexes of $A^{e}$-modules $\tau_{\bullet} : C_{\bullet}(A) \rightarrow K_{\bullet}(A)$ 
and $\bar{\tau}_{\bullet} : \bar{C}_{\bullet}(A) \rightarrow K_{\bullet}(A)$, giving quasi-isomorphisms, by the semisimplicity assumption on $k^{e}$. 

Let us start with $\tau_{1} b_{2} = d_{2} \tau_{2}$. 
\begin{itemize}
\item If we apply $d_{2} \tau_{2}$ to a normalized element of the form $1|a|b|1$, we see that $d_{2} \tau_{2}(1|a|b|1) = 0$. 
On the other hand,  
\begin{align*} 
   \tau_{1} b_{2} (1|a|b|1) &= \tau_{1}(a|b|1-1|ab|1+1|a|b) 
   \\ 
   &= a b_{(1)}|\bar{b}_{(2)}|b_{(3)} - (ab)_{(1)}|\overline{(ab)}_{(2)}|(ab)_{(3)} + a_{(1)}|\bar{a}_{(2)}|a_{(3)} b,     
\end{align*}
which trivially vanishes by Fact \ref{fac:cuenta division}, so $\tau_{1} b_{2} = d_{2} \tau_{2}$ for the elements $1|a|b|1$, with $\deg(a) + \deg(b) < N$. 

\item If we apply $d_{2} \tau_{2}$ to a relation decomposition of the form $1|a_{i}|b_{i}|1$ (with $r = \sum_{i \in I} a_{i} \otimes b_{i}$), 
we see that $d_{2} \tau_{2}(1|a_{i}|b_{i}|1) = d_{2}(1|r|1) = r_{(1)}|\bar{r}_{(2)}|r_{(3)}$. 
Analogously,  
\begin{align*} 
    \tau_{1} b_{2} (1|a_{i}|b_{i}|1) &= \tau_{1}(a_{i}|b_{i}|1-1|a_{i}b_{i}|1+1|a_{i}|b_{i}) 
    \\
       &= a_{i} b_{i,(1)}|\bar{b}_{i,(2)}|b_{i,(3)} + a_{i,(1)}|\bar{a}_{i,(2)}|a_{i,(3)} b_{i},     
\end{align*} 
since $1|a_{i}b_{i}|1 = 1|r|1 \in C_{1}(A)$, which vanishes. 
Now, using that 
\[     r_{(1)}|\bar{r}_{(2)}|r_{(3)} = a_{i} b_{i,(1)}|\bar{b}_{i,(2)}|b_{i,(3)} + a_{i,(1)}|\bar{a}_{i,(2)}|a_{i,(3)} b_{i},     \]
we conclude that $\tau_{1} b_{2} = d_{2} \tau_{2}$ where we have defined it. 
\end{itemize}

Let us now prove that $\tau_{2} b_{3} = d_{3} \tau_{3}$. 
\begin{itemize}
\item If we apply $d_{3} \tau_{3}$ to a normalized element of the form $1|a|b|c|1$, we see that $d_{3} \tau_{3}(1|a|b|c|1) = 0$. 
On the other hand,  
\[     \tau_{2} b_{3} (1|a|b|c|1) = \tau_{2}(a|b|c|1-1|ab|c|1+1|a|bc|1-1|a|b|c),     \]
which trivially vanishes by \eqref{eq:tau2,1}, so $\tau_{2} b_{3} = d_{3} \tau_{3}$ for the elements of the form $1|a|b|c|1$, 
with $\deg(a) + \deg(b) +\deg(c) < N$. 

\item If we apply $d_{3} \tau_{3}$ to a relation decomposition of the form $1|a_{i}|b_{i}|c_{i}|1$, we see that $d_{3} \tau_{3}(1|a_{i}|b_{i}|c_{i}|1) = 0$. 
Also, 
\[     \tau_{2} b_{3} (1|a_{i}|b_{i}|c_{i}|1) = \tau_{2}(a_{i}|b_{i}|c_{i}|1-1|a_{i}b_{i}|c_{i}|1+1|a_{i}|b_{i}c_{i}|1-1|a_{i}|b_{i}|c_{i}) 
       = 0,     \]
where we have used \eqref{eq:tau2,1}, since $\deg(a_{i}), \deg(c_{i}) > 0$.

\item If we apply $d_{3} \tau_{3}$ to a double relation decomposition of the form $1|a_{i}|b_{i}|c_{i}|1$ 
(with $w = \sum_{i \in I} a_{i} \otimes b_{i} \otimes c_{i}$), 
we see that $d_{3} \tau_{3}(1|a_{i}|b_{i}|c_{i}|1) = d_{3}(1|w|1) = a_{i}|b_{i} c_{i}|1 - 1|a_{i} b_{i}|c_{i}$. 
We note that, if $w = r_{i} u_{i} = v_{i} s_{i}$, then $a_{i}|b_{i} c_{i}|1 = v_{i}|s_{i}|1$ and $1|a_{i} b_{i}|c_{i} = 1|r_{i}|u_{i}$, 
since, for the first identity, each member is a decomposition of $w$ in two terms such that first one has degree $1$, 
and the argument for the second identity is analogous.
On the other hand,  
\begin{align*}
   \tau_{2} b_{3} (1|a_{i}|b_{i}|c_{i}|1) &= \tau_{2}(a_{i}|b_{i}|c_{i}|1-1|a_{i} b_{i}|c_{i}|1+1|a_{i}|b_{i} c_{i}|1-1|a_{i}|b_{i}|c_{i}) 
   \\
    &= a_{i}|b_{i} c_{i}|1 - 1|a_{i} b_{i}|c_{i},     
\end{align*} 
where we have used that $1|a_{i} b_{i}|c_{i}|1 = 1|r_{i}|u_{i}|1$ and $1|a_{i}|b_{i} c_{i}|1 = 1|v_{i}|s_{i}|1$ vanish in $C_{2}(A)$, 
and \eqref{eq:tau2,2} for the other two terms. 
We conclude that $\tau_{2} b_{3} = d_{3} \tau_{3}$, where we have defined it. 
\end{itemize}

Finally, there exist homotopies $s_{\bullet}$ and $t_{\bullet}$ for $\bullet \in \NN_{0}$, which are morphisms of $A$-bimodules, on the complexes $\bar{C}_{\bullet}(A)$ 
and $K_{\bullet}(A)$, respectively, such that 
\begin{equation}
\label{eq:sn}
     \bar{b}_{n+1} s_{n} + s_{n-1} \bar{b}_{n} = 1_{\bar{C}_{n}(A)} - \bar{\sigma}_{n} \bar{\tau}_{n},     
\end{equation}
and 
\begin{equation}
\label{eq:tn}
     d_{n+1} t_{n} + t_{n-1} d_{n} = 1_{K_{n}(A)} - \bar{\tau}_{n} \bar{\sigma}_{n},     
\end{equation}
hold for $n \in \NN_{0}$, respectively. 
We set
\begin{align}
   &s_{-1} = 0, 
   \label{eq:s-1}\addtocounter{equation}{1}\tag{S$_{-1}$} 
   \\
   &s_{0} = 0, 
   \label{eq:s0}\addtocounter{equation}{1}\tag{S$_{0}$} 
   \\
   &s_{1}(1|a|1) = - 1|a_{(1)}|\bar{a}_{(2)}|a_{(3)}, \text{if $\deg(a) < N$,}
   \label{eq:s1}\addtocounter{equation}{1}\tag{S$_{1}$} 
   \\
   &\begin{multlined}
   s_{2}(1|a_{i}|b_{i}|1) = 1|a_{i}|b_{i,(1)}|\bar{b}_{i,(2)}|b_{i,(3)}, 
   \\
   \text{if $1|a_{i}|b_{i}|1$ is normalized or a relation decomposition,}
   \end{multlined}
   \label{eq:s2}\addtocounter{equation}{1}\tag{S$_{2}$}
   \\
   &\begin{multlined}
   s_{3}(1|a_{i}|b_{i}|c_{i}|1) = - 1|a_{i}|b_{i}|c_{i,(1)}|\bar{c}_{i,(2)}|c_{i,(3)}, 
   \\
   \text{if $1|a_{i}|b_{i}|c_{i}|1$ is normalized, a rel. decomp. or a double rel. decomp.}
   \end{multlined}
   \label{eq:s3,3}\addtocounter{equation}{1}\tag{S$_{3}$}
\end{align}
Note that in fact $s_{3}(1|a_{i}|b_{i}|c_{i}|1)$ vanishes if $1|a_{i}|b_{i}|c_{i}|1$ is a double relation decomposition.

As previously indicated, the semisimplicity hypothesis on $k^{e}$ tells us that these maps can be extended to complete morphisms of $A^{e}$-modules 
$s_{\bullet} : \bar{C}_{\bullet}(A) \rightarrow \bar{C}_{\bullet+1}(A)$ satisfying the identity \eqref{eq:sn}. 
We see clearly that $1_{\bar{C}_{0}(A)} - \bar{\sigma}_{0} \bar{\tau}_{0} = \bar{b}_{1} s_{0} + s_{-1} \bar{b}_{0}$ holds. 
Moreover, \eqref{eq:sn} for $n = 1$ is also verified, since, for $\deg(a) < N$,  
\[     (1_{\bar{C}_{1}(A)} - \bar{\sigma}_{1} \bar{\tau}_{1})(1|a|1) = 1|a|1 - a_{(1)}|\bar{a}_{(2)}|a_{(3)},     \] 
and
\begin{align*}
   (\bar{b}_{2} s_{1} + s_{0} \bar{b}_{1})(1|a|1) &= \bar{b}_{2} (
                                                - 1|a_{(1)}|\bar{a}_{(2)}|a_{(3)}) 
   \\                                 &= 
                                         - a_{(1)}|\bar{a}_{(2)}|a_{(3)} + 1|a_{(1)} \bar{a}_{(2)}|a_{(3)} - 1|a_{(1)}|\bar{a}_{(2)} a_{(3)}   
   \\
                                      &= 
                                         - a_{(1)}|\bar{a}_{(2)}|a_{(3)} + 1|a|1 
                                       = 1|a|1 - a_{(1)}|\bar{a}_{(2)}|a_{(3)}, 
\end{align*}
where we have used Fact \ref{fac:cuenta telescopica} in the penultimate equality. 

We shall now check that \eqref{eq:sn} holds for $n = 2$ when applied to a normalized $1|a|b|1$, 
or to a relation decomposition $1|a_{i}|b_{i}|1$. 
\begin{itemize}
\item If we apply \eqref{eq:sn} for $n = 2$ to a normalized element of the form $1|a|b|1$, 
we see that $(1_{\bar{C}_{2}(A)} - \bar{\sigma}_{2} \bar{\tau}_{2})(1|a|b|1) = 1|a|b|1$, since $\bar{\tau}_{2}(1|a|b|1)=0$. 
On the other hand,  
\begin{align*}
     s_{1} \bar{b}_{2} (1|a|b|1) &= s_{1}(a|b|1-1|ab|1+1|a|b) 
     \\                    &= 
                              - a|b_{(1)}|\bar{b}_{(2)}|b_{(3)} 
                             + 1|(ab)_{(1)}|\overline{(ab)}_{(2)}|(ab)_{(3)} 
                             - 1|a_{(1)}|\bar{a}_{(2)}|a_{(3)} b
     \\
                           &= 
                              - a|b_{(1)}|\bar{b}_{(2)}|b_{(3)} + 1|a b_{(1)}|\bar{b}_{(2)}|b_{(3)},     
\end{align*}
where we have used Fact \ref{fac:cuenta division}. 
Besides, 
\begin{align*}
     \bar{b}_{3} s_{2} (1|a|b|1) &= \bar{b}_{3}(1|a|b_{(1)}|\bar{b}_{(2)}|b_{(3)} 
                              ) 
     \\
                           &= a|b_{(1)}|\bar{b}_{(2)}|b_{(3)} - 1|a b_{(1)}|\bar{b}_{(2)}|b_{(3)} + 1|a|b_{(1)}\bar{b}_{(2)}|b_{(3)}
     \\
                           &- 1|a|b_{(1)}|\bar{b}_{(2)}b_{(3)} 
     \\                       
                            &= a|b_{(1)}|\bar{b}_{(2)}|b_{(3)} - 1|a b_{(1)}|\bar{b}_{(2)}|b_{(3)} + 1|a|b|1,
\end{align*}
where we have used Fact \ref{fac:cuenta telescopica} in the last equality. 
By adding the previous computations, we see that \eqref{eq:sn} holds for $n = 2$ when applied to $1|a|b|1$ if $\deg(a) + \deg(b) < N$. 

\item If we apply \eqref{eq:sn} for $n = 2$ to a relation decomposition $1|a_{i}|b_{i}|1$ (with $r = \sum_{i \in I} a_{i} \otimes b_{i}$), 
we see that 
\begin{multline*}    
     (1_{\bar{C}_{2}(A)} - \bar{\sigma}_{2} \bar{\tau}_{2})(1|a_{i}|b_{i}|1) = 1|a_{i}|b_{i}|1 - \bar{\sigma}_{2}(1|r|1) 
\\  
      = 1|a_{i}|b_{i}|1 - 1|a_{i} b_{i,(1)}|\bar{b}_{i,(2)}|b_{i,(3)} - 1|a_{i,(1)}|\bar{a}_{i,(2)}|a_{i,(3)} b_{i}.     
\end{multline*}
On the other hand,  
\begin{align*} 
     s_{1} \bar{b}_{2} (1|a_{i}|b_{i}|1) &= s_{1}(a_{i}|b_{i}|1+1|a_{i}|b_{i}) 
     \\
     &= 
       - a_{i}|b_{i,(1)}|\bar{b}_{i,(2)}|b_{i,(3)} - 1|a_{i,(1)}|\bar{a}_{i,(2)}|a_{i,(3)} b_{i},     
\end{align*} 
where we have used that the elements $1|a_{i}b_{i}|1 = 1|r|1$ 
in $\bar{C}_{1}(A)$ vanish.
Also, 
\begin{align*}
     \bar{b}_{3} s_{2} (1|a_{i}|b_{i}|1) &= \bar{b}_{3}(1|a_{i}|b_{i,(1)}|\bar{b}_{i,(2)}|b_{i,(3)}
                                    ) 
     \\
                            &= a_{i}|b_{i,(1)}|\bar{b}_{i,(2)}|b_{i,(3)} - 1|a_{i} b_{i,(1)}|\bar{b}_{i,(2)}|b_{i,(3)} 
     \\
                            &+ 1|a_{i}|b_{i,(1)}\bar{b}_{i,(2)}|b_{i,(3)}
                            - 1|a_{i}|b_{i,(1)}|\bar{b}_{i,(2)}b_{i,(3)} 
     \\                       
                            &= a_{i}|b_{i,(1)}|\bar{b}_{i,(2)}|b_{i,(3)} - 1|a_{i} b_{i,(1)}|\bar{b}_{i,(2)}|b_{i,(3)} + 1|a_{i}|b_{i}|1,
\end{align*}
where we have used Fact \ref{fac:cuenta telescopica} in the last equality. 
Adding both computations, we see that \eqref{eq:sn} holds for $n = 2$ when applied to $1|a_{i}|b_{i}|1$, with $r = \sum_{i \in I} a_{i} \otimes b_{i} \in R$, $\deg(a_{i})$, $\deg(b_{i}) > 0$.
\end{itemize}

We shall now check that \eqref{eq:sn} holds for $n = 3$ when applied to a normalized $1|a|b|c|1$, 
to a relation decomposition or to a double relation decomposition $1|a_{i}|b_{i}|c_{i}|1$. 
\begin{itemize}
\item If we apply \eqref{eq:sn} for $n = 3$ to a normalized element or to a relation decomposition of the form $1|a_{i}|b_{i}|c_{i}|1$, 
we see that $(1_{\bar{C}_{3}(A)} - \bar{\sigma}_{3} \bar{\tau}_{3})(1|a_{i}|b_{i}|c_{i}|1) = 1|a_{i}|b_{i}|c_{i}|1$, since $\bar{\tau}_{3}(1|a_{i}|b_{i}|c_{i}|1)=0$. 
On the other hand,  
\begin{align*}
     s_{2} \bar{b}_{3} (1|a_{i}|b_{i}|c_{i}|1) &= s_{2}(a_{i}|b_{i}|c_{i}|1 - 1|a_{i} b_{i}|c_{i}|1 + 1|a_{i}|b_{i} c_{i}|1 - 1|a_{i}|b_{i}|c_{i}) 
     \\                    
                                         &= a_{i}|b_{i}|c_{i,(1)}|\bar{c}_{i,(2)}|c_{i,(3)} 
                                          - 1|a_{i} b_{i}|c_{i,(1)}|\bar{c}_{i,(2)}|c_{i,(3)} 
     \\
                                         &+ 1|a_{i}|(b_{i}c_{i})_{(1)}|\overline{(b_{i}c_{i})}_{(2)}|(b_{i}c_{i})_{(3)} 
                                          - 1|a_{i}|b_{i,(1)}|\bar{b}_{i,(2)}|b_{i,(3)} c_{i} 
     \\
                                         &= a_{i}|b_{i}|c_{i,(1)}|\bar{c}_{i,(2)}|c_{i,(3)} 
                                          - 1|a_{i} b_{i}|c_{i,(1)}|\bar{c}_{i,(2)}|c_{i,(3)} 
     \\                                     
                                         &+  1|a_{i}|b_{i,(1)}|\bar{b}_{i,(2)}|b_{i,(3)} c_{i} + 1|a_{i}|b_{i}c_{i,(1)}|\bar{c}_{i,(2)}|c_{i,(3)}  
     \\
                                         &- 1|a_{i}|b_{i,(1)}|\bar{b}_{i,(2)}|b_{i,(3)} c_{i}
     \\
                                         &= a_{i}|b_{i}|c_{i,(1)}|\bar{c}_{i,(2)}|c_{i,(3)} 
                                          - 1|a_{i} b_{i}|c_{i,(1)}|\bar{c}_{i,(2)}|c_{i,(3)} 
     \\
                                         &+ 1|a_{i}|b_{i}c_{i,(1)}|\bar{c}_{i,(2)}|c_{i,(3)},                                               
\end{align*}
where we have used Fact \ref{fac:cuenta division}. 
Besides, 
\begin{align*}
     \bar{b}_{4} s_{3} (1|a_{i}|b_{i}|c_{i}|1) &= \bar{b}_{4}(
                                                  -1|a_{i}|b_{i}|c_{i,(1)}|\bar{c}_{i,(2)}|c_{i,(3)}) 
     \\
                           &= 
                           - a_{i}|b_{i}|c_{i,(1)}|\bar{c}_{i,(2)}|c_{i,(3)} + 1|a_{i} b_{i}|c_{i,(1)}|\bar{c}_{i,(2)}|c_{i,(3)}
     \\          
                           &- 1|a_{i}|b_{i} c_{i,(1)}|\bar{c}_{i,(2)}|c_{i,(3)} + 1|a_{i}|b_{i}|c_{i}|1, 
\end{align*}
where we have used Fact \ref{fac:cuenta telescopica} in the last equality. 
The addition of these computations tells us that \eqref{eq:sn} holds for $n = 2$ when applied to a normalized element or a relation decomposition.

\item If we apply \eqref{eq:sn} for $n = 3$ to a double relation decomposition $1|a_{i}|b_{i}|c_{i}|1$ 
(with $w = \sum_{i \in I} a_{i} \otimes b_{i} \otimes c_{i}$, $w=r_{i}u_{i}=v_{i}s_{i}$), 
we see that 
\begin{align*}
     (1_{\bar{C}_{3}(A)} - \bar{\sigma}_{3} \bar{\tau}_{3})(1|a_{i}|b_{i}|c_{i}|1) &= 1|a_{i}|b_{i}|c_{i}|1 - \bar{\sigma}_{3}(1|w|1) 
     \\
            &= 1|a_{i}|b_{i}|c_{i}|1 - 1|a_{i}|(b_{i}c_{i})_{(1)}| \overline{(b_{i}c_{i})}_{(2)}| (b_{i}c_{i})_{(3)} 
     \\     
            &= 1|a_{i}|b_{i}|c_{i}|1 - 1|a_{i}|b_{i,(1)}| \bar{b}_{i,(2)}| b_{i,(3)} c_{i} 
     \\
            &- 1|a_{i}|b_{i} c_{i,(1)}| \bar{c}_{i,(2)}| c_{i,(3)}
     \\
            &= - 1|a_{i}|b_{i,(1)}| \bar{b}_{i,(2)}| b_{i,(3)} c_{i},
\end{align*}                                                                   
where we have used that $1|a_{i}|b_{i} c_{i,(1)}| \bar{c}_{i,(2)}| c_{i,(3)} = 1|a_{i}|b_{i}|c_{i}|1$, for the degree of $c_{i}$ is one. 
On the other hand,  
\begin{align*}
     s_{2} \bar{b}_{3} (1|a_{i}|b_{i}|c_{i}|1) &= s_{2}(a_{i}|b_{i}|c_{i}|1 - 1|a_{i}|b_{i}|c_{i}) 
     \\  
     &= 
        - 1|a_{i}|b_{i,(1)}|\bar{b}_{i,(2)}|b_{i,(3)} c_{i},     
\end{align*}
where we have used that $1|a_{i} b_{i}|c_{i}|1 = 1|r_{i}|u_{i}|1$ and $1|a_{i}|b_{i} c_{i}|1 = 1|v_{i}|s_{i}|1$ vanish in $\bar{C}_{2}(A)$, 
and $a_{i}|b_{i}|c_{i,(1)}|\bar{c}_{i,(2)}|c_{i,(3)}$ vanish in $\bar{C}_{3}(A)$. 
Also, 
\begin{align*}
     \bar{b}_{4} s_{3} (1|a_{i}|b_{i}|c_{i}|1) &= 0,
\end{align*}
using that $s_{3}$ vanishes on a double relation decomposition. 
We conclude thus that \eqref{eq:sn} holds for $n = 3$ when applied to a double relation decomposition. 
\end{itemize}

The comparison morphism $\bar{\sigma}_{\bullet} : K_{\bullet}(A) \rightarrow \bar{C}_{\bullet}(A)$ is an injection, for $\bullet = 0, \dots, 3$.  
Furthermore, the morphisms $\bar{\sigma}_{\bullet}$ and $\bar{\tau}_{\bullet}$ satisfy the following result.
\begin{lemma}
\label{lem:tau sigma} 
According to the previous definition of the comparison morphisms $\bar{\sigma}_{\bullet} : K_{\bullet}(A) \rightarrow \bar{C}_{\bullet}(A)$ 
and $\bar{\tau}_{\bullet} : \bar{C}_{\bullet}(A) \rightarrow K_{\bullet}(A)$ for $\bullet = 0, \dots, 3$, we see that 
$\bar{\tau}_{\bullet} \bar{\sigma}_{\bullet} = 1_{K_{\bullet}(A)}$.
\end{lemma}
\noindent\textbf{Proof.} 
It is clear that $\bar{\tau}_{0} \bar{\sigma}_{0} = 1_{A} \otimes 1_{A}$ and $\bar{\tau}_{1} \bar{\sigma}_{1} = 1_{A} \otimes 1_{V} \otimes 1_{A}$. 
Moreover, 
\begin{align*}
     \bar{\tau}_{2} \bar{\sigma}_{2} (1|r|1) &= \bar{\tau}_{2}(1|r_{(1)}|\bar{r}_{(2)}|r_{(3)})   
     \\
                                 &= \bar{\tau}_{2}(\underset{\deg(r_{(3)})>0}{\underbrace{1|r_{(1)}|\bar{r}_{(2)}|r_{(3)}}} + 1|r_{(1)}|\bar{r}_{(2)}|1)
     \\                            
                                 &= 0 + \bar{\tau}_{2}(1|r_{(1)}|\bar{r}_{(2)}|1) = 1|r|1.     
\end{align*}

Finally, 
\begin{align*}
     \bar{\tau}_{3} \bar{\sigma}_{3} (1|w|1) &= \bar{\tau}_{3}(1|v_{i}|s_{i,(1)}|\bar{s}_{i,(2)}|s_{i,(3)})   
     \\
                                 &= \bar{\tau}_{3}(\underset{\deg(s_{i,(3)})>1}{\underbrace{1|v_{i}|s_{i,(1)}|\bar{s}_{i,(2)}|s_{i,(3)}}} 
                                 + \underset{\deg(s_{i,(3)})=1}{\underbrace{1|v_{i}|s_{i,(1)}|\bar{s}_{i,(2)}|s_{i,(3)}}}
                                 + 1|v_{i}|s_{i,(1)}|\bar{s}_{i,(2)}|1)
     \\                            
                                 &= 0 + \bar{\tau}_{3}(1|\bar{r}_{i,(1)}|r_{i,(2)}|\bar{r}_{i,(3)}|u_{i}) + \bar{\tau}_{3}(1|v_{i}|s_{i,(1)}|\bar{s}_{i,(2)}|1) 
                                 = 1|w|1,   
\end{align*}
where $w = v_{i} s_{i} = r_{i} u_{i} \in R_{N+1}$ and we have used that 
$1|v_{i}|s_{i,(1)}|\bar{s}_{i,(2)}|\bar{s}_{i,(3)} = 1|\bar{r}_{i,(1)}|r_{i,(2)}|\bar{r}_{i,(3)}|u_{i}$. 
\qed

\begin{remark}
\label{rem:t}
The Lemma implies that we may choose $t_{\bullet} = 0$, for $\bullet = -1, \dots, 2$.
\end{remark}

We have thus obtained a partial description of both projective resolutions and their comparison in lower degrees 
that can be depicted as follows:
\begin{small}
\[
\xymatrix@C-13pt
{
\dots
\ar[r]
&
A \otimes R_{N+1} \otimes A
\ar[r]^-{d_{3}}
\ar@<+.7ex>[d]^{\bar{\sigma}_{3}}
\ar@/^/[l]^-{t_{3}}
&
A \otimes R \otimes A
\ar[r]^-{d_{2}}
\ar@<+.7ex>[d]^{\bar{\sigma}_{2}}
\ar@/^/[l]^-{t_{2}}
&
A \otimes V \otimes A
\ar[r]^-{d_{1}}
\ar@<+.7ex>[d]^{\bar{\sigma}_{1}}
\ar@/^/[l]^-{t_{1}}
&
A \otimes A
\ar[r]^-{d_{0}}
\ar@<+.7ex>[d]^{\bar{\sigma}_{0}}
\ar@/^/[l]^-{t_{0}}
&
A
\ar[r]
\ar@{=}[d]
\ar@/^/[l]^-{t_{-1}}
&
0
\\
\dots
\ar[r]
&
A \otimes (A/k)^{\otimes 3} \otimes A
\ar[r]^-{\bar{b}_{3}}
\ar@<+.7ex>[u]^{\bar{\tau}_{3}}
\ar@/^/[l]^-{s_{3}}
&
A \otimes (A/k)^{\otimes 2} \otimes A
\ar[r]^-{\bar{b}_{2}}
\ar@<+.7ex>[u]^{\bar{\tau}_{2}}
\ar@/^/[l]^-{s_{2}}
&
A \otimes (A/k) \otimes A
\ar[r]^-{\bar{b}_{1}}
\ar@<+.7ex>[u]^{\bar{\tau}_{1}}
\ar@/^/[l]^-{s_{1}}
&
A \otimes A
\ar[r]^-{b_{0}}
\ar@<+.7ex>[u]^{\bar{\tau}_{0}}
\ar@/^/[l]^-{s_{0}}
&
A
\ar[r]
\ar@/^/[l]^-{s_{-1}}
&
0
}
\]
\end{small}

\section{PBW-deformations and deformations \`a la Gerstenhaber of homogeneous algebras}
\label{sec:pbwger}

In this section we shall briefly recall the definitions of PBW-deformations and of the (graded) deformations \emph{\`a la Gerstenhaber}, 
which we will usually just call deformations. 
We shall also establish a link between both concepts. 

\subsection{PBW-deformations of homogeneous algebras}
\label{subsec:pbw}

We start recalling the definition of a PBW-deformation, and we mainly follow \cite{BG06}. 
We first recall that a \emph{filtered $k$-algebra} $B$ is a $k$-algebra provided with an increasing sequence $\{ F^{\bullet}B \}_{\bullet \in \NN_{0}}$ of $k$-subbimodules of $B$ 
such that $F^{m}B.F^{n}B \subseteq F^{m+n}B$, for all $m,n \in \NN_{0}$, and $1_{B} \in F^{0}B$. 
As usual, such filtrations may also be seen to be indexed over $\ZZ$, where the negatively indexed terms vanish. 
Given a $k$-bimodule $V$, the tensor algebra $TV$ has a filtration $\{F^{\bullet}\}_{\bullet \in \NN_{0}}$ defined by $F^{i} = \oplus_{j=0}^{i} V^{\otimes j}$. 
Now, given $P \subset F^{N}$, we shall consider the algebra $U = TV/\cl{P}$, with the filtration $\{F^{\bullet}U\}_{\bullet \in \NN_{0}}$ induced by the filtration of the tensor algebra, 
\textit{i.e.} $F^{\bullet}U = \pi(F^{\bullet})$, where $\pi$ denotes the canonical projection from $TV$ to $U$. 
Of course, $\pi$ is a morphism of filtered algebras.
The filtration can be described more concretely as follows: 
if $\cl{P}^{i} = F^{i} \cap \cl{P}$, then $F^{i}U = F^{i}/\cl{P}^{i}$, for $i \in \NN_{0}$. 
If $\pi_{i} : TV \rightarrow V^{\otimes i}$ is the canonical projection, let us denote $R = \pi_{N}(P)$ and define the $N$-homogeneous algebra 
$A = TV/\cl{R}$. 

\begin{remark}
We remark the standard fact that, even though the ideal $\cl{P}$ coincides with $\sum_{i,j \geq 0} V^{\otimes i} P V^{\otimes j}$, 
$\cl{P}^{n}$ may be strictly bigger than the sum $\sum_{i + j \leq n-N} V^{\otimes i} P V^{\otimes j}$, which in particular vanishes if $n < N$. 
As we shall see below, the PBW property will be the exact condition in order to avoid this phenomenon. 
\end{remark}

We shall now consider the associated graded algebra $\mathrm{gr}(U)$ to the previous filtration. 
First, we state the following direct results.
\begin{lemma}
\label{lema:filgrad}
If $S = \oplus_{i \in \NN_{0}} S_{i}$ is an $\NN_{0}$-graded $k$-algebra and $F^{\bullet}S$ is the filtration induced by the grading of $S$, 
\textit{i.e.} $F^{i}S = \oplus_{j=0}^{i} S_{j}$, then there exists a canonical isomorphism $\iota : S \rightarrow \mathrm{gr}(S)$ of $\NN_{0}$-graded $k$-algebras, 
such that the restriction $\iota|_{S_{i}} : S_{i} \rightarrow F^{i}S/F^{i-1}S$ is the canonical $k^{e}$-linear isomorphism.
\end{lemma}
\noindent\textbf{Proof.} 
Easy.
\qed


Since $\pi : TV \rightarrow U$ is a morphism of filtered algebras, it induces a morphism of graded algebras $\mathrm{gr}(\pi) : \mathrm{gr}(TV) \rightarrow \mathrm{gr}(U)$. 
Moreover, the filtration of $U$ is induced by the filtration of $TV$, so $\mathrm{gr}(\pi)$ is surjective. 
On the other hand, since the filtration of $TV$ comes from a grading on the tensor algebra, we see that there exists a canonical isomorphism $\iota : TV \simeq \mathrm{gr}(TV)$, 
by Lemma \ref{lema:filgrad}.  
So we may consider the surjective morphism of graded $k$-algebras given by the composition $\mathrm{gr}(\pi) \circ \iota : TV \rightarrow \mathrm{gr}(U)$, which we shall call $\Pi$. 
It is easy to see that $\Pi(R) = 0$, since $\iota(R) = P/F^{N-1}$. 
Hence $\Pi$ induces a surjective morphism of graded $k$-algebras $p : A \rightarrow \mathrm{gr}(U)$. 
We say that $U$ satisfies the \emph{PBW property} or that $U$ is a \emph{PBW-deformation} of $A$ if $p$ is an isomorphism. 

\begin{remark}
If $k = F$ is a field and $V$ is a finite dimensional vector space over $k$, $p$ is an isomorphism if and only if there exists an isomorphism of graded 
$k$-algebras $A \simeq \mathrm{gr}(U)$. 
This is proved as follows. 
One direction is obvious. 
Let us assume that there exists an isomorphism of graded $k$-algebras $A \simeq \mathrm{gr}(U)$. 
Since $V$ is finite dimensional, $A$ (and \textit{a fortiori} $\mathrm{gr}(U)$) is evidently locally finite dimensional, \textit{i.e.} each homogeneous component is finite dimensional. 
Hence, since $p$ is surjective, each restriction of $p$ to a homogeneous component is surjective, so bijective. 
Therefore, $p$ is an isomorphism.
\end{remark}

The morphism $p$ can be more concretely described as follows. 
We just need to consider its restriction to $A_{i}$ ($i \in \NN_{0}$). 
First, we see that 
\[     \mathrm{gr}(U)_{i} = F^{i}U/F^{i-1}U = \frac{(F^{i}/(F^{i} \cap \cl{P}))}{(F^{i-1}/(F^{i-1} \cap \cl{P}))} 
                          \simeq F^{i}/((F^{i} \cap \cl{P})+F^{i-1}),      \]
where we have used that $F^{i-1} \cap \cl{P} = (F^{i} \cap \cl{P}) \cap F^{i-1}$ and the Second and Third Isomorphism Theorems. 
Then, $p|_{A_{i}}$ is induced by the $k^{e}$-linear morphism $V^{\otimes i} \rightarrow F^{i}/((F^{i} \cap \cl{P})+F^{i-1}) \simeq \mathrm{gr}(U)_{i}$ 
given by the composition of the canonical injection $V^{\otimes i} \rightarrow F^{i}$ and the canonical projection $F^{i} \rightarrow F^{i}/((F^{i} \cap \cl{P})+F^{i-1})$. 
So it is easy to see that $U$ satisfies the PBW property if and only if 
\[     (\cl{P} \cap F^{n}) + F^{n-1} = (\cl{R} \cap F^{n}) + F^{n-1},     \]
for all $n \in \NN_{0}$, which is equivalent to
\[     \cl{P} \cap F^{n} \subset (\cl{R} \cap F^{n}) + F^{n-1},     \]
for all $n \in \NN_{0}$.

As noticed by R. Berger and V. Ginzburg (see \cite{BG06}, Prop. 3.2), the filtered algebra $U$ satisfies the PBW property if and only if 
$\cl{P}^{n} = \sum_{i + j \leq n-N} V^{\otimes i} P V^{\otimes j}$, for all $n \in \NN_{0}$ 
(in fact, it is sufficient to prove the equality for $n \geq N-1$). 
Moreover, if we denote $J_{n} = \sum_{i + j \leq n-N} V^{\otimes i} P V^{\otimes j}$, for $n \in \NN_{0}$, Prop. 3.3. in \cite{BG06} 
states that $U$ satisfies the PBW property if and only if $J_{n} \cap F^{n-1} = J_{n-1}$, for all $n \in \NN_{0}$ (or just $n \geq N$). 
The identity $J_{N} \cap F^{N-1} = J_{N-1}$ is simply 
\begin{equation}
\label{eq:pbw1}
     P \cap F^{N-1} = 0,     
\end{equation}
whereas $J_{N+1} \cap F^{N} = J_{N}$ is easily equivalent to
\begin{equation}
\label{eq:pbw2}
     (V \otimes P + P \otimes V) \cap F^{N} \subset P.     
\end{equation}

From now on we shall suppose that identity \eqref{eq:pbw1} holds, which implies that the map $\pi_{N} : F^{N} \rightarrow V^{\otimes N}$ gives an isomorphism between $P$ and $R = \pi_{N}(P)$. 
Then there exists a $k^{e}$-linear map $\varphi : R \rightarrow F^{N-1}$ such that $\mathrm{id}  - \varphi$ is the inverse of $\pi_{N}|_{P}$, 
\textit{i.e.} $P = \{ r - \varphi(r) : r \in R \}$. 
We further write, $\varphi = \sum_{j = 0}^{N-1} \varphi_{j}$, where $\varphi_{j} : R \rightarrow V^{\otimes j}$ is the composition of $\varphi$ with the canonical morphism $F^{N-1} \rightarrow  V^{\otimes j}$. 
Then it is easy to see that identity \eqref{eq:pbw2} is equivalent to (see \cite{BG06}, Prop. 3.5)
\[     (\varphi \otimes 1_{V} - 1_{V} \otimes \varphi) (R_{N+1}) \subset P,     \]
or equivalently (see \cite{BG06}, Prop. 3.6)
\begin{align}
     &(\varphi_{N-1} \otimes 1_{V} - 1_{V} \otimes \varphi_{N-1})(R_{N+1}) \subset R,
     \label{eq:pbwfi1}
     \\
     &\varphi_{0} \circ (\varphi_{N-1} \otimes 1_{V} - 1_{V} \otimes \varphi_{N-1})(R_{N+1}) = 0,
     \label{eq:pbwfi3}
     \\
     &(\varphi_{j} \circ (\varphi_{N-1} \otimes 1_{V} - 1_{V} \otimes \varphi_{N-1}) + (\varphi_{j-1} \otimes 1_{V} - 1_{V} \otimes \varphi_{j-1}))(R_{N+1}) = 0, 
     \label{eq:pbwfi2}
\end{align}
for $0 < j< N$. 

\begin{definition}
Given a filtered algebra $U=TV/\cl{P}$, where $P \subset F^{N}$, such that \eqref{eq:pbw1}, \eqref{eq:pbwfi1}, \eqref{eq:pbwfi2} and \eqref{eq:pbwfi3} hold, we say that 
$U$ is a \emph{weak PBW-deformation} of $A = TV/\cl{R}$, where $R = \pi_{N}(P)$. 
We remark that each weak PBW-deformation $U$ of $A$ is provided with a surjective morphism of graded algebras $p : A \rightarrow \mathrm{gr}(U)$. 
Given two weak deformations $U$ and $U'$ of $A$ (with induced morphisms $p$ and $p'$, respectively), they are called \emph{equivalent} if there exists an isomorphism of filtered algebras 
$g : U \rightarrow U'$ such that $\mathrm{gr}(g) \circ p = p'$. 
\end{definition}

It is immediate to see that a PBW-deformation is a weak PBW-deformation. 
Using considerations of deformation theory \emph{\`a la Gerstenhaber} we shall provide another proof of the converse as a consequence of Theorem \ref{teo:pbw} when $A$ is an $N$-homogeneous algebra satisfying that $\Tor_{3}^{A}(k,k)$ is concentrated in degree $N+1$ (\textit{cf.}~\cite{BG06}, Thm. 3.4). 

\subsection{Graded deformations in the sense of Gerstenhaber}
\label{subsec:ger}

Let us now briefly recall the definition of a graded deformation and some results that we shall use in the sequel. 
Most of what we will present is implicit in the work of M. Gerstenhaber (see \cite{Ge64}, Sec. 1.2-1.5), and it is explained in more detail 
by A. Braverman and D. Gaitsgory in \cite{BGa96}. 
We would like to remark, however, that we are working over a not necessarily commutative ring $k$ and this situation needs more sophisticated tools 
(\textit{cf.} \cite{Gor08}, Sec. 2, where the author is dealing with the case $k = \CC[G]$ or more generally $k$ a separable $\CC$-algebra, 
in order to assure that the Hochschild cohomology over $k$ coincides with that over $\CC$). 

In what follows, we consider $k[t]$ as an $\NN_{0}$-graded $k$-algebra such that $\deg(t) = 1$ and $t$ is central. 
Given a $k$-bimodule $V$, we shall denote by $V[t]$ the $k[t]$-bimodule with elements $\sum_{i \in I} v_{i} t^{i}$, for $v_{i} \in V$ and 
finite subsets $I \subseteq \NN_{0}$, provided with the action 
$c t^{m} (\sum_{i \in I} v_{i} t^{i}) c' t^{n} = \sum_{i \in I} c.v_{i}.c' t^{i+m+n}$, for $c, c' \in k$. 

If $A$ denotes an $\NN_{0}$-graded associative $k$-algebra and $i \in \NN$, an \emph{$i$-th level graded deformation} of $A$ means a graded $k[t]/(t^{i+1})$-algebra structure on the $k[t]/(t^{i+1})$-bimodule $A_{i} = A[t]/(t^{i+1})$ such that the identity $A_{i}/t.A_{i} \simeq A$ is an isomorphism of graded algebras. 
By a \emph{(polynomial) graded deformation} of $A$ we mean a graded $k[t]$-algebra structure on the $k[t]$-bimodule $A_{t} = A[t]$ 
such that the identity $A_{t}/t.A_{t} \simeq A$ is an isomorphism of graded algebras. 
In the previous definitions, we are always using the obvious graded $k^{e}$-linear map $A \rightarrow A_{i} = A[t]/(t^{i+1})$ (resp., $A \rightarrow A_{t} = A[t]$) given by $a \mapsto a$. 
We will usually denote the product of $A_{i}$ (resp. $A_{t}$) by $\times^{i}$ (resp. $\times$), which can be written as
\[     a \times^{i} b = a b + \sum_{h = 1}^{i} \psi_{h}(a,b) t^{h} \hskip 0.5cm \text{($a \times b = a b + \sum_{h \in \NN} \psi_{h}(a,b) t^{h}$)}.     \] 
Since $\psi_{h}$ has degree $-h$, we must note that the sum for $\times$ is finite for any pair of homogeneous elements $a$ and $b$ in $A$. 
We remark that $\psi_{h}$ may be considered as an element of $\mathcal{H}om_{A^{e}}(C_{2}(A),A)$, and that the associativity of 
$\times$ is equivalent to 
\begin{align}
   d\psi_{1} &= 0,
  \label{eq:Gerdef1}
   \\
   - d\psi_{j+1}(a,b,c) &= \sum_{i=1}^{j} (\psi_{i}(a,\psi_{j+1-i}(b,c)) - \psi_{i}(\psi_{j+1-i}(a,b),c)), \text{for $j \in \NN$.}
   \label{eq:Gerdef2}
\end{align}
The right member of the last equation is usually denoted by $\mathrm{sq}(\psi_{1},\dots,\psi_{j})(a,b,c)$. 
We note that there exists a trivial polynomial deformation of $A$ given by the trivial product on $A[t]$, \textit{i.e.} such that 
$(a t^{m}) \times_{0} (a' t^{n}) = (a.a') t^{m+n}$, for $a, a' \in A$ and $m,n \in \NN_{0}$. 

Given a filtered algebra $B$ with filtration $\{ F^{\bullet}B \}_{\bullet \in \NN_{0}}$, we recall that the \emph{Rees algebra} associated to it is the graded $k[t]$-algebra 
\[     R(B) = \{ \sum_{i \in I} b_{i} t^{i} : \hskip 0.6mm \text{$I$ is finite and $b_{i} \in F^{i}B$} \},    \]
which is considered as a subalgebra of $B[t]$ provided with the trivial product $\times_{0}$ (in this case $B$ is concentrated in degree zero). 
We remark that the underlying graded algebra structure of $R(B)$ is $\oplus_{\bullet \in \NN_{0}} F^{\bullet}B$ with the product induced by that of $B$. 
It is easy to see that $R(B)/\cl{t-\lambda} \simeq B$, for $\lambda \in k^{\times} \cap \Z(k)$ such that its action on $B$ is central 
(\textit{i.e.} $\lambda b = b \lambda$, for all $b \in B$) and $R(B)/\cl{t} \simeq \mathrm{gr}(B)$ (\textit{cf.}~\cite{CG97}, Cor. 2.3.8, 
whose proof can be applied also to this case).
Moreover, $R(\place)$ defines a functor from the category of filtered $k$-algebras to the category of graded $k[t]$-algebras. 

We would like to make some remarks about the algebra $A_{t}/\cl{t-1}$ (or about the algebra $A_{t}/\cl{t-\lambda}$, 
with $\lambda \in k^{\times} \cap \Z(k)$ such that its action on $A$ is central, to which the following arguments also apply). 
There exists a $k^{e}$-linear map $A \rightarrow A_{t}/\cl{t-1}$ given by the composition of the canonical injection $A \rightarrow A_{t}$ and the projection $A_{t} \rightarrow A_{t}/\cl{t-1}$. 
We consider the filtration on $A_{t}/\cl{t-1}$ induced by the filtration of $A$ under the previous map. 
We remark that the filtration of $A_{t}/\cl{t-1}$ induced by the filtration of $A_{t}$ that comes from the grading is trivial. 
We shall see that the associated graded algebra of $A_{t}/\cl{t-1}$ is isomorphic to $A$ as graded algebras. 
This is proved as follows. 
We consider the $k^{e}$-linear map $\rho' : A_{t} \rightarrow A$ (not an algebra map) given by $\sum_{j=0}^{m} a_{j} t^{j} \mapsto \sum_{j=0}^{m} a_{j}$. 
It is trivially verified that $\rho'$ respects the filtrations coming from the gradings, and that $\rho'((t-1) b)=0$, for any $b \in A_{t}$, so it induces a morphism of filtered $k$-bimodules 
$\rho^{*} : A_{t}/\cl{t-1} \rightarrow A$, which is obviously surjective and injective. 
Its inverse is just the aforementioned map $A \rightarrow A_{t}/\cl{t-1}$. 
Hence it induces an isomorphism of graded $k$-bimodules $\rho : \mathrm{gr}(A_{t}/\cl{t-1}) \rightarrow \mathrm{gr}(A) \simeq A$ 
(the last isomorphism by Lemma \ref{lema:filgrad}). 
Now, if we denote the product of $A_{t}/\cl{t-1}$ by $\times_{1}$ and taking into account that the product of two elements 
$a, b \in A$ in $A_{t}/\cl{t-1}$ is given by 
\[     a \times_{1} b = a.b + \sum_{j \in \NN} \psi_{j}(a,b),     \]
and the degree of $\psi_{j}$ is $-j$, we see that $\rho$ is an algebra morphism, so an isomorphism of graded algebras. 
This further implies that, if $A$ is generated by the image of the $k^{e}$-linear map $V \rightarrow A$, then $A_{t}/\cl{t-1}$ is generated by the image of the composition of $V \rightarrow A$ and 
the $k^{e}$-linear map $A \rightarrow A_{t}/\cl{t-1}$. 

The importance of the algebra $A_{t}/\cl{t-1}$ is explained in the following proposition. 
\begin{proposition}
\label{prop:fibra}
Let $A$ be a graded algebra and let $A_{t}$ be a graded deformation of $A$. 
Then, there exists a canonical isomorphism of graded $k[t]$-algebras $R(A_{t}/\cl{t-1}) \rightarrow A_{t}$, such that the induced morphism 
$A \simeq R(A_{t}/\cl{t-1})/\cl{t} \rightarrow A_{t}/\cl{t} \simeq A$ is the identity, where $A \simeq R(A_{t}/\cl{t-1})/\cl{t}$ is given by the composition of the inverse of $\rho$ and the canonical 
isomorphism of graded algebras $\mathrm{gr}(A_{t}/\cl{t-1}) \simeq R(A_{t}/\cl{t-1})/\cl{t}$. 
\end{proposition}
\noindent\textbf{Proof.}
We first remark that we are going to identify $A_{t}/\cl{t-1}$ with $A$ as $k^{e}$-bimodules under the previous filtered $k^{e}$-linear isomorphism 
$A \rightarrow A_{t}/\cl{t-1}$, 
so we will denote the elements of this last algebra by elements $a$ of $A$. 
We will continue denoting the product of $A_{t}/\cl{t-1}$ by $\times_{1}$ and we remark that, if $a, b \in A$ are of degree $i$ and $j$ respectively, then
\[     a \times_{1} b = a b + \sum_{h = 1}^{i+j} \psi_{h}(a,b),     \]
since $\psi_{h}$ has degree $-h$. 
Under this identification, the homogeneous $i$-th component of $R(A_{t}/\cl{t-1})$ is just $\oplus_{j=0}^{i} A_{j}$. 
So the elements of $R(A_{t}/\cl{t-1})$ are sums of elements of the form $a_{i} t^{j}$ for $j \geq i \geq 0$, with $a_{i} \in A_{i}$. 

We recall that the elements of $A_{t}$ are sums of elements of the form $a_{i} t^{j}$ for $i, j \geq 0$, for $a_{i} \in A_{i}$. 

We now define the map $\mathrm{com} : R(A_{t}/\cl{t-1}) \rightarrow A_{t}$ given by the linear extension of 
\[     a_{i} t^{j} \mapsto a_{i} t^{j-i}.     \]
It is trivially verified that $\mathrm{com}$ is an isomorphism of graded $k[t]$-bimodules. 
In order to prove that it is a morphism of algebras, we only need to show that 
\[     \mathrm{com}((a_{i} t^{j}) \times_{1} (b_{l} t^{m})) = \mathrm{com}(a_{i} t^{j}) \times \mathrm{com}(b_{l} t^{m}),     \]
for $a_{i}, b_{l} \in A$ of degree $i$ and $l$ respectively, and $i \leq j$ and $l \leq m$.
The left member is given by 
\begin{align*}
     \mathrm{com}((a_{i} t^{j}) \times_{1} (b_{l} t^{m})) &= \mathrm{com} ((a_{i} \times_{1} b_{l}) t^{j+m}) 
     = \mathrm{com} ((a_{i} b_{l} + \sum_{h = 1}^{i+l} \psi_{h}(a_{i},b_{l})) t^{j+m}) 
     \\
     &= \mathrm{com} (a_{i} b_{l} t^{j+m}) + \sum_{h = 1}^{i+l} \mathrm{com}(\psi_{h}(a_{i},b_{l}) t^{j+m}) 
     \\
     &= a_{i} b_{l} t^{j+m-i-l} + \sum_{h = 1}^{i+l} \psi_{h}(a_{i},b_{l}) t^{j+m+h-i-l},
\end{align*}
where we have used that $\psi_{h}$ has degree $-h$. 
On the other hand, the right member is given by 
\begin{align*}
     \mathrm{com}(a_{i} t^{j}) \times \mathrm{com}(b_{l} t^{m}) &= (a_{i} t^{j-i}) \times (b_{l} t^{m-l})
     \\
     &= (a_{i} \times b_{l}) t^{j+m-i-l} = ((a_{i} b_{l} + \sum_{h = 1}^{i+l} \psi_{h}(a_{i},b_{l}) t^{h}) t^{j+m-i-l} 
     \\
     &= a_{i} b_{l} t^{j+m-i-l} + \sum_{h = 1}^{i+l} \psi_{h}(a_{i},b_{l}) t^{j+m+h-i-l},
\end{align*}
where we have again used that $\psi_{h}$ has degree $-h$. 
Hence, $\mathrm{com}$ is an isomorphism of graded $k[t]$-algebras. 
It is clear that the induced map $A \simeq R(A_{t}/\cl{t-1})/\cl{t} \rightarrow A_{t}/\cl{t} \simeq A$ is the identity. 
\qed

Let $E(A)$ denote the groupoid of all graded deformations of $A$ where the morphisms are by definition isomorphisms of graded $k[t]$-algebras, 
such that the induced morphism of $A \simeq A_{t}/t.A_{t}$ is the identity. 
Analogously, let $E_{i}(A)$ denote the groupoid of all $i$-th level graded deformations of $A$. 
Given $i \in \NN$, we denote by $F_{i}$ the functor from $E(A)$ to $E_{i}(A)$ given by reduction modulo $t^{i+1}$, \textit{i.e.} if $A_{t}$ is a graded deformation of $A$, then 
$F_{i}(A_{t}) = A_{t}/t^{i+1}A_{t}$ and the definition on morphisms is the obvious one.
Moreover, given $i < j$ natural numbers, we denote by $F_{i < j}$ the functor from $E_{j}(A)$ to $E_{i}(A)$ given by reduction modulo $t^{i+1}$, so if $A_{j}$ denotes a $j$-th level deformation, 
$F_{i < j}(A_{j}) = A_{j}/t^{i+1}A_{j}$ and for the morphisms it is obvious. 

The following lemma is trivial (\textit{cf.}~\cite{BGa96}, Lemma 1.3).
\begin{lemma}
The collection of functors $F_{i}$ define an equivalence between the category $E(A)$ and the inverse limit of the categories $E_{i}(A)$ with respect to the functors $F_{i<j}$.
\end{lemma}

Given an $i$-th level deformation $A_{i}$ of $A$, a \emph{continuation to an $(i+1)$-th level deformation of $A_{i}$} is an $(i+1)$-th level deformation $A_{i+1}$ 
of $A$ such that $F_{i<i+1}(A_{i+1}) = A_{i}$. 
Given two continuations $A_{i+1}$ and $A'_{i+1}$ of $A_{i}$ to an $(i+1)$-th level deformation, a morphism $f$ from $A_{i+1}$ to $A'_{i+1}$ is a morphism in $E_{i+1}(A)$ such that 
$F_{i<i+1}(f) = 1_{A_{i}}$. 
The following proposition is also immediate (\textit{cf.}~\cite{BGa96}, Prop. 1.5, or \cite{Gor08}, 2.6).
\begin{proposition}
\label{prop:basic}
\begin{itemize}
\item[(a)] The set of isomorphism classes of objects of $E_{1}(A)$ can be canonically identified with $HH^{2}(A)_{-1}$.

\item[(b)] Let $A_{i}$ be an object of $E_{i}(A)$. 
Then the obstruction for its continuation to an $(i + 1)$-th level deformation lies in $HH^{3}(A)_{-i-1}$.

\item[(c)] Let $A_{i}$ be as in (b). 
Then the set of isomorphism classes of continuations of $A_{i}$ to an $(i + 1)$-th level deformation has a natural structure of an $HH^{2}(A)_{-i-1}$-homogeneous space.
\end{itemize}
\end{proposition}

Finally we state the following proposition, which is analogous to Prop. 3.7 of \cite{BGa96}. 
\begin{proposition}
\label{prop:koszul}
Let $A$ be an $N$-homogeneous algebra that is Koszul. 
Then,
\begin{itemize}
\item[(i)] the functors $F_{i-1<i} : E_{i}(A) \rightarrow E_{i-1}(A)$ are injective on isomorphism classes of
objects for $i > N$,
\item[(ii)] the functors $F_{i-1<i}$ are surjective on isomorphism classes of objects for $i > N+1$.
\end{itemize}
\end{proposition}
\noindent\textbf{Proof.}
It is easy to see from the bimodule Koszul complex that $HH^{2}_{-i}(A)$ vanishes for
$i > N$. 
Hence Proposition \ref{prop:basic}(c) implies (i). 
Analogously, $HH^{3}_{-i}(A)$ vanishes for $i > N+1$,
so Proposition \ref{prop:basic}(b) implies (ii).
\qed

In what follows, we shall only consider graded deformations such that the unit of the original $k$-algebra is also a unit of the deformed algebra. 
This is equivalent to ask that the $2$-cochains $\psi_{j}$ actually belong to $\mathcal{H}om_{A^{e}}(\bar{C}_{\bullet}(A),A)$. 
We shall say that such graded deformations \emph{preserve the unit}. 
Since any graded deformation is equivalent to another one preserving the unit (because the complexes $\mathcal{H}om_{A^{e}}(\bar{C}_{\bullet}(A),A)$ and 
$\mathcal{H}om_{A^{e}}(C_{\bullet}(A),A)$ are quasi-isomorphic, and the equations \eqref{eq:Gerdef1} and \eqref{eq:Gerdef2} for both complexes 
are preserved under the corresponding quasi-isomorphism), there is would be no loss of generality in restricting to such a situation. 

\subsection{From deformations to weak PBW-deformations}
\label{subsec:gertopbw}

After having recalled the basic definitions and results we shall provide a link between both concepts. 
From now on, we stress the fact that we assume that that the graded $k$-algebra $A$ is of the form $TV/\cl{R}$, 
where $V$ is considered to be concentrated in degree $1$ and $R \subseteq V^{\otimes N}$, 
for $N \geq 2$ satisfying the assumption that $\Tor_{3}^{A}(k,k)$ is concentrated in degree $N + 1$, so we may use the considerations of Section \ref{sec:gen}.
First, we set some notation: if $\psi : W \rightarrow A$ denotes a $k^{e}$-linear map, we shall denote $\psi^{\sim} : A \otimes W \otimes A \rightarrow A$, its unique $A^{e}$-linear extension. 
A $k^{e}$-linear map $\psi : (A/k)^{\otimes 2} \rightarrow A$ is called \emph{normalized} if it vanishes on normalized elements, 
\textit{i.e.} if $\psi(a \otimes b) = 0$, whenever $a,b \in I_{+} \simeq A/k$ are homogeneous elements satisfying that $\deg(a) + \deg(b) < N$. 
Further, we say that $\psi$ is \emph{extranormalized} if it is normalized and if it vanishes on relation decompositions. 
We shall also say in this case that $\psi^{\sim}$ is normalized or extranormalized. 

Let us now consider a graded deformation of $A$ given by 
\[     a \times b = a.b + \sum_{i \geq 1} \psi_{i}(a,b) t^{i},     \] 
where $\psi_{i} : (A/k)^{\otimes 2} \rightarrow A$ 
are normalized $k^{e}$-linear maps 
(see Lemma \ref{lema:sim}). 
We remark that $\psi_{i}$ is a morphism of degree $-i$, and we have that the associativity of $\times$ is equivalent to \eqref{eq:Gerdef1} and \eqref{eq:Gerdef2}. 

We set $\varphi_{N-j}^{\sim} = \psi_{j}^{\sim} \circ \bar{\sigma}_{2}$, for $j = 1,\dots, N$. 
We point out that the grading implies that $\psi_{j}^{\sim} \circ \bar{\sigma}_{2}$ vanishes for $j > N$. 
We shall see that, once we assume that \eqref{eq:Gerdef1} and \eqref{eq:Gerdef2} hold, the induced morphisms $\varphi_{N-j} : R \rightarrow F^{N-1}$ satisfy the weak PBW property expressed in 
\eqref{eq:pbwfi1}, \eqref{eq:pbwfi2} and \eqref{eq:pbwfi3}. 

Let us first prove \eqref{eq:pbwfi1}. 
We shall state a simple fact that we shall use in the sequel. 
\begin{fact}
\label{fac:difkos}
Given $\varphi : R \rightarrow A$ a $k^{e}$-linear map, then $d\varphi : R_{N+1} \rightarrow A$ satisfies that
\[     d\varphi = (1_{V} \otimes \varphi - \varphi \otimes 1_{V}),     \]
where we remark that the map $(1_{V} \otimes \varphi - \varphi \otimes 1_{V})$ is defined from $R_{N+1}$ to $A$. 
\end{fact}
\noindent\textbf{Proof.}
If $w = r_{i} u_{i} = v_{i} s_{i} \in R_{N+1}$, then 
\begin{align*}
      d\varphi(w) &= d\varphi^{\sim}(1|w|1) = (\varphi^{\sim} \circ d_{3}) (1|w|1) = \varphi^{\sim} (v_{i}|r_{i}|1-1|s_{i}|u_{i}) 
      \\
                  &= v_{i} \varphi(r_{i})- \varphi(s_{i})u_{i} = (1_{V} \otimes \varphi - \varphi \otimes 1_{V})(w).     
\end{align*}
\qed

We note that \eqref{eq:Gerdef1} implies the equality 
\[     d\varphi_{N-1}^{\sim} = d(\psi_{1}^{\sim} \circ \bar{\sigma}_{2}) = \psi_{1}^{\sim} \circ \bar{\sigma}_{2} \circ d_{3} 
                             = \psi_{1}^{\sim} \circ \bar{b}_{3} \circ \bar{\sigma}_{3} = d\psi_{1}^{\sim} \circ \bar{\sigma}_{3} = 0.     \]
Hence, Fact \ref{fac:difkos} implies that $(1_{V} \otimes \varphi_{N-1} - \varphi_{N-1} \otimes 1_{V})(w) = 0$. 
So, when we consider $(1_{V} \otimes \varphi_{N-1} - \varphi_{N-1} \otimes 1_{V})$ as a map from $R_{N+1}$ to $V^{\otimes N}$, we get \eqref{eq:pbwfi1}. 
More precisely, we see that equation \eqref{eq:Gerdef1} composed with $\bar{\sigma}_{3}$ is equivalent to equation \eqref{eq:pbwfi1}.

Let us now prove \eqref{eq:pbwfi2} and \eqref{eq:pbwfi3}. 
We shall need the following simple fact.
\begin{fact}
\label{fac:conor}
If $\psi : (A/k)^{\otimes 2} \rightarrow A$ is a normalized cocycle and $\gamma \in V^{\otimes N}$, then it holds that 
$\psi(\bar{\gamma}_{(1)},\gamma_{(2)}) = \psi(\gamma_{(1)},\bar{\gamma}_{(2)})$. 
More generally, let us suppose that $\psi_{1}, \dots, \psi_{j} : (A/k)^{\otimes 2} \rightarrow A$ are normalized $k^{e}$-linear maps such that $-d\psi_{j} = \mathrm{sq}(\psi_{1},\dots,\psi_{j-1})$. 
Hence, if $\gamma \in V^{\otimes N}$, we see that $\psi_{j}(\bar{\gamma}_{(1)},\gamma_{(2)}) = \psi_{j}(\gamma_{(1)},\bar{\gamma}_{(2)})$. 
\end{fact}
\noindent\textbf{Proof.}
It is obvious that the second statement generalizes the first one, but we give a detailed proof of both. 
In the first case we note that, since $\psi$ is a cocycle, 
\begin{align*}
     0 &= d\psi(\bar{\gamma}_{(1)},\gamma_{(2)},\bar{\gamma}_{(3)}) 
         = \bar{\gamma}_{(1)} \psi(\gamma_{(2)},\bar{\gamma}_{(3)}) - \psi(\bar{\gamma}_{(1)}\gamma_{(2)},\bar{\gamma}_{(3)}) 
     \\
         &\phantom{= d\psi(\bar{\gamma}_{(1)},\gamma_{(2)},\bar{\gamma}_{(3)}) l}
            + \psi(\bar{\gamma}_{(1)},\gamma_{(2)}\bar{\gamma}_{(3)}) - \psi(\bar{\gamma}_{(1)},\gamma_{(2)})\bar{\gamma}_{(3)} 
     \\
        &= - \psi(\gamma_{(1)},\bar{\gamma}_{(2)}) + \psi(\bar{\gamma}_{(1)},\gamma_{(2)}),
\end{align*}
where we have used the normalization of $\psi$ in the third equality.

Now we prove the second statement. 
On the one hand, just as before we have that 
\begin{align*}
     d\psi_{j}(\bar{\gamma}_{(1)},\gamma_{(2)},\bar{\gamma}_{(3)}) 
     &= \bar{\gamma}_{(1)} \psi_{j}(\gamma_{(2)},\bar{\gamma}_{(3)}) - \psi_{j}(\bar{\gamma}_{(1)}\gamma_{(2)},\bar{\gamma}_{(3)}) 
     \\
     &+ \psi_{j}(\bar{\gamma}_{(1)},\gamma_{(2)}\bar{\gamma}_{(3)}) - \psi_{j}(\bar{\gamma}_{(1)},\gamma_{(2)})\bar{\gamma}_{(3)} 
     \\
        &= - \psi_{j}(\gamma_{(1)},\bar{\gamma}_{(2)}) + \psi_{j}(\bar{\gamma}_{(1)},\gamma_{(2)}),
\end{align*}
where we have used the normalization of $\psi_{j}$. 
On the other hand, 
\begin{multline*}
     \mathrm{sq}(\psi_{1},\dots,\psi_{j-1})(\bar{\gamma}_{(1)},\gamma_{(2)},\bar{\gamma}_{(3)}) 
      \\ = \sum_{i=1}^{j-1} (\psi_{i}(\bar{\gamma}_{(1)},\psi_{j-i}(\gamma_{(2)},\bar{\gamma}_{(3)})) - \psi_{i}(\psi_{j-i}(\bar{\gamma}_{(1)},\gamma_{(2)}),\bar{\gamma}_{(3)}))
      = 0,     
\end{multline*}
since $\psi_{1},\dots,\psi_{j-1}$ are normalized. 
Hence our statement follows. 
\qed

Using Fact \ref{fac:difkos} for $\varphi_{N-j}$ ($0 < j < N$) we see that 
\begin{equation}
\label{eq:mitad1}     
(d\psi_{j+1}^{\sim} \circ \bar{\sigma}_{3}) (1|w|1) =  (1_{V} \otimes \varphi_{N-j-1} - \varphi_{N-j-1} \otimes 1_{V})(w).     
\end{equation}
Besides, $(\mathrm{sq}(\psi_{1},\dots,\psi_{j})^{\sim} \circ \bar{\sigma}_{3}) (1|w|1)$ is equal to 
\begin{align*}
   \mathrm{sq}(\psi_{1},\dots,&\psi_{j})^{\sim} (1|v_{i}|s_{i,(1)}|\bar{s}_{i,(2)}|s_{i,(3)}) 
   \\
   &= \sum_{i=1}^{j} (\psi_{i}(v_{i},\psi_{j+1-i}(s_{i,(1)},\bar{s}_{i,(2)})) - \psi_{i}(\psi_{j+1-i}(v_{i},s_{i,(1)}),\bar{s}_{i,(2)}))s_{i,(3)} 
   \\
   &= \sum_{i=1}^{j} (\psi_{i}(v_{i},\psi_{j+1-i}(s_{i,(1)},\bar{s}_{i,(2)})) - \psi_{i}(\psi_{j+1-i}(v_{i},s_{i,(1)}),\bar{s}_{i,(2)}))
   \\
   &= \sum_{i=1}^{j} (\psi_{i}(v_{i},\psi_{j+1-i}(s_{i,(1)},\bar{s}_{i,(2)})) - \psi_{i}(\psi_{j+1-i}(\bar{r}_{i,(1)},r_{i,(2)}),u_{i}))
   \\
   &= \psi_{j}(v_{i},\psi_{1}(s_{i,(1)},\bar{s}_{i,(2)})) - \psi_{j}(\psi_{1}(\bar{r}_{i,(1)},r_{i,(2)}),u_{i}),
\end{align*}
where we have used in the third equality the normalization of $\psi_{1}, \dots, \psi_{j}$. 
In the last equality we have used the following simple fact: $\psi_{j+1-i}$ has degree $-(j+1-i) = i-j-1$, so 
$\psi_{j+1-i}$ applied to an element of degree $N$ gives an element of degree $N+i-j-1$. 
In consequence, the elements  $v_{i} \otimes \psi_{j+1-i}(s_{i,(1)},\bar{s}_{i,(2)})$ and $\psi_{j+1-i}(\bar{r}_{i,(1)},r_{i,(2)}) \otimes u_{i}$ have degree $N+j-i \leq N$ 
(for $i \leq j$) and the degree is exactly $N$ if and only if $i=j$. 
The normalization of $\psi_{j}$ forces only to consider the terms with degree $N$ and the last equality follows. 

Now, using Fact \ref{fac:conor} we see that $\psi_{1}(\bar{r}_{i,(1)},r_{i,(2)}) = \psi_{1}(r_{i,(1)},\bar{r}_{i,(2)})$. 
Furthermore, by its very definition, $\varphi_{N-1}(r) = \psi_{1}(r_{(1)}, \bar{r}_{(2)})$, for any $r \in R$. 
Therefore, 
\[     (\mathrm{sq}(\psi_{1},\dots,\psi_{j})^{\sim} \circ \bar{\sigma}_{3}) (1|w|1) 
         = \psi_{j}(v_{i},\varphi_{N-1}(s_{i})) - \psi_{j}(\varphi_{N-1}(r_{i}),u_{i}).     \]
We need to compare the previous expression with $\varphi_{N-j} \circ (1_{V} \circ \varphi_{N-1} - \varphi_{N-1} \otimes 1_{V})(w)$, as we shall proceed to do. 
We first note that $\varphi_{N-j} \circ (1_{V} \otimes \varphi_{N-1} - \varphi_{N-1} \otimes 1_{V})(w)$ is given by 
\begin{align*}
   \varphi_{N-j}(v_{i} \otimes &\varphi_{N-1}(s_{i})) - \varphi_{N-j}(\varphi_{N-1}(r_{i}) \otimes u_{i})
   \\
   &= (\psi_{j}^{\sim} \circ \bar{\sigma}_{2}) (1|\overset{\phantom{r} \in R}{\overbrace{v_{i} \otimes \underset{\alpha_{i}}{\underbrace{\varphi_{N-1}(s_{i})}} 
   - \underset{\beta_{i}}{\underbrace{\varphi_{N-1}(r_{i})}} \otimes u_{i}}}|1)
   \\
   &= \psi_{j}^{\sim} (1|v_{i}\alpha_{i,(1)}|\bar{\alpha}_{i,(2)}|\alpha_{i,(3)} - 1|\beta_{i,(1)}|\bar{\beta}_{i,(2)}|\beta_{i,(3)}u_{i} - 1|\beta_{i}|u_{i}|1)
   \\
   &= \psi_{j}(v_{i}\alpha_{i,(1)},\bar{\alpha}_{i,(2)})\alpha_{i,(3)} - \psi_{j}(\beta_{i,(1)},\bar{\beta}_{i,(2)})\beta_{i,(3)}u_{i} - \psi_{j}(\beta_{i},u_{i})
   \\
   &= \psi_{j}(v_{i}\alpha_{i,(1)},\bar{\alpha}_{i,(2)}) - \psi_{j}(\beta_{i},u_{i})
    = \psi_{j}(v_{i},\alpha_{i}) - \psi_{j}(\beta_{i},u_{i})
   \\
   &= \psi_{j}(v_{i},\varphi_{N-1}(s_{i})) - \psi_{j}(\varphi_{N-1}(r_{i}),u_{i}),
\end{align*}
where we have considered in the third member that $v_{i} \otimes \varphi_{N-1}(s_{i}) - \varphi_{N-1}(r_{i}) \otimes u_{i}$ is an element of $R$, we have used the normalization of  
$\psi_{j}$ in the fourth equality, and Fact \ref{fac:conor} in the first term of the penultimate member. 
This implies that 
\begin{equation}
\label{eq:mitad2}
     (\mathrm{sq}(\psi_{1},\dots,\psi_{j})^{\sim} \circ \bar{\sigma}_{3}) (1|w|1) = \varphi_{N-j} \circ (1_{V} \otimes \varphi_{N-1} - \varphi_{N-1} \otimes 1_{V})(w).   
\end{equation}

Now, from equations \eqref{eq:mitad1} and \eqref{eq:mitad2} we see that equation \eqref{eq:Gerdef2} composed with $\bar{\sigma}_{3}$ (for $j = 1, \dots, N-1$) is exactly \eqref{eq:pbwfi2}. 
The case of identity \eqref{eq:Gerdef2} for $j=N$ leads exactly to equation \eqref{eq:pbwfi3}, since 
$d\psi_{N+1}^{\sim} \circ \bar{\sigma}_{3} = 0$ by degree reasons. 

We have thus proved the following result.
\begin{proposition}
\label{prop:gertopbw}
Let $\times$ be a graded deformation of $A$ given by 
\[     a \times b = a.b + \sum_{i \geq 1} \psi_{i}(a,b) t^{i},     \] 
where $\psi_{i} : (A/k)^{\otimes 2} \rightarrow A$ (resp. $\psi_{i}^{\sim} : \bar{C}_{2}(A) \rightarrow A$) are normalized $k^{e}$-linear maps 
(resp. $A^{e}$-linear maps). 
We define $\varphi_{N-j}^{\sim} = \psi_{j}^{\sim} \circ \bar{\sigma}_{2}$, for $j = 1,\dots, N$. 
Then, the induced morphisms $\varphi_{N-j} : R \rightarrow F^{N-1}$ satisfy the weak PBW property expressed in 
\eqref{eq:pbwfi1}, \eqref{eq:pbwfi2} and \eqref{eq:pbwfi3}.
\end{proposition}

Even though we have considered graded deformations given by normalized maps, the following lemma shows that the assumption is in fact unnecessary.
\begin{lemma}
\label{lema:sim}
Let $A$ be an $N$-homogeneous algebra such that $\Tor_{3}^{A}(k,k)$ is concentrated in degree $N+1$. 
Then, any graded deformation of $A$ is equivalent to a deformation preserving the unit given by normalized maps. 
\end{lemma}
\noindent\textbf{Proof.}
Let $A_{t}$ be a graded deformation of $A$ given by a collection of maps $\{ \psi_{j} : (A/k)^{\otimes 2} \rightarrow A \}_{j \in \NN}$ 
(each of degree $- j$) which define a product $\times$. 
We only need to show that there exists another graded deformation $A_{t}'$ given by maps $\{ \psi'_{j} : (A/k)^{\otimes 2} \rightarrow A \}_{j \in \NN}$ 
such that the first $N$ maps $\psi'_{1}, \dots, \psi'_{N}$ are normalized, for the $\{ \psi'_{j} \}_{j > N}$ are always automatically normalized 
by degree reasons. 

Using the main property \eqref{eq:sn} of the homotopy $s_{\bullet}$, we see that 
$\psi_{1}^{\sim} - \psi_{1}^{\sim} \circ \bar{\sigma}_{2} \circ \bar{\tau}_{2} = d(\psi_{1}^{\sim} \circ s_{1})$. 
It is clear that $\psi_{1}^{\sim} \circ \bar{\sigma}_{2} \circ \bar{\tau}_{2}$ is normalized. 
Let us define $\alpha_{1} : \bar{A} \rightarrow A$ the map induced by $\psi_{1}^{\sim} \circ s_{1}$. 
So we see that $\mathrm{exp}(t \alpha_{1})$ gives an equivalence from the algebra $(A_{t},\times)$ to 
another deformation $(A_{t}^{1},\times_{1})$ of $A$ given by the maps $\{ \psi^{1}_{\bullet} : A^{\otimes 2} \rightarrow A \}_{\bullet \in \NN}$ 
such that $(\psi^{1}_{1})^{\sim} = \psi_{1}^{\sim} \circ \bar{\sigma}_{2} \circ \bar{\tau}_{2}$ is normalized. 

We proceed now by (finite) induction. 
Let $i \leq N$ and let us suppose that $A_{t}$ is equivalent to a deformation $A_{t}^{i}$ given by a collection of maps 
$\{ \psi^{i}_{\bullet} : (A/k)^{\otimes 2} \rightarrow A \}_{\bullet \in \NN}$ such that $\psi^{i}_{\bullet}$ are normalized for $\bullet \leq i$. 
Using again \eqref{eq:sn}, we conclude that
\[     (\psi^{i}_{i+1})^{\sim} - 
       ((\psi^{i}_{i+1})^{\sim} \circ \bar{\sigma}_{2} \circ \bar{\tau}_{2} - \mathrm{sq}(\psi^{i}_{1}, \dots , \psi^{i}_{i})^{\sim} \circ s_{2}) 
       = d((\psi^{i}_{i+1})^{\sim} \circ s_{1}).     \] 
Also, it is clear that 
$(\psi^{i}_{i+1})^{\sim} \circ \bar{\sigma}_{2} \circ \bar{\tau}_{2} - \mathrm{sq}(\psi^{i}_{1}, \dots , \psi^{i}_{i})^{\sim} \circ s_{2}$ 
is normalized (each summand is obviously so). 
Let us define $\alpha_{i+1} : \bar{A} \rightarrow A$ the map induced by $(\psi^{i}_{i+1 })^{\sim} \circ s_{1}$. 
This tells us that $\mathrm{exp}(t^{i+1} \alpha_{i+1})$ gives an equivalence from the algebra $(A^{i}_{t},\times_{i})$ to 
another deformation $(A_{t}^{i+1},\times_{i+1})$ of $A$ given by the maps $\{ \psi^{i+1}_{\bullet} : (A/k)^{\otimes 2} \rightarrow A \}_{\bullet \in \NN}$ 
such that $\psi^{i+1}_{j} = \psi^{i}_{j}$ for $j \leq i$ 
and the map $(\psi^{i+1}_{i+1})^{\sim} = 
(\psi^{i}_{i+1})^{\sim} \circ \bar{\sigma}_{2} \circ \bar{\tau}_{2} - \mathrm{sq}(\psi^{i}_{1}, \dots , \psi^{i}_{i})^{\sim} \circ s_{2}$ is normalized. 
\qed

\subsection{From weak PBW-deformations to deformations}
\label{subsec:pbwtoger}

Now, we shall give an inverse construction, which is a little more complicated. 
Let us suppose that $\varphi = \sum_{j=0}^{N-1} \varphi_{j}$, $\varphi_{j} : R \rightarrow V^{\otimes j}$, gives a weak PBW-deformation of an $N$-homogeneous algebra $A$ satisfying that $\Tor_{3}^{A}(k,k)$ is concentrated in degree $N+1$. 
We shall construct a (possible infinite) sequence of normalized maps $\{\psi_{j}\}_{j \in \NN}$, where $\psi_{j} : (A/k)^{\otimes 2} \rightarrow A$, 
such that  
\[     a \times b = a.b + \sum_{i \geq 1} \psi_{i}(a,b) t^{i}     \] 
is a graded deformation of $A$. 

First, we define $(\psi'_{j})^{\sim} = \varphi_{N-j}^{\sim} \circ \bar{\tau}_{2}$ for $j = 1, \dots, N$, and zero otherwise. 
We note that $\psi'_{j}$ is a normalized homogeneous morphism of degree $-j$ for $j = 1, \dots, N$. 

We note that $\psi'_{1}$ is a cocycle since 
\[     d(\psi'_{1})^{\sim} = (\psi'_{1})^{\sim} \circ \bar{b}_{3} = \varphi_{N-1}^{\sim} \circ \bar{\tau}_{2} \circ \bar{b}_{3} 
                        = \varphi_{N-1}^{\sim} \circ d_{3} \circ \bar{\tau}_{3} = d\varphi_{N-1}^{\sim} \circ \bar{\tau}_{3},     \]
and Fact \ref{fac:difkos} tells us that the evaluation of the last expression at $1|w|1$, for $w \in R_{N+1}$, is equal to the element of $A$ given by 
$(1_{V} \otimes \varphi_{N-1} - \varphi_{N-1} \otimes 1_{V})(w)$. 
The identity \eqref{eq:pbwfi1} says that this element 
vanishes, so $d(\psi'_{1})^{\sim} = 0$. 
We define $\psi_{1} = \psi'_{1}$. 

We will proceed recursively on $j \in \NN$. 
Let us suppose that we have defined $\psi_{1}, \dots, \psi_{j}$ such that $\psi_{i}^{\sim} - (\psi'_{i})^{\sim}$ is an extranormalized 
$k^{e}$-linear map of degree $-i$ for all $i = 1, \dots, j$, and that 
\[        - d\psi_{i+1} = \mathrm{sq}(\psi_{1},\dots,\psi_{i})     \]
holds for all $i=0, \dots, j-1$. 
We will denote $\eta_{i}^{\sim} = \psi_{i}^{\sim} - (\psi'_{i})^{\sim}$.

We shall now prove that 
\begin{equation}
\label{eq:Gdefs3}
        - d(\psi'_{j+1})^{\sim} \circ \bar{\sigma}_{3} = \mathrm{sq}(\psi_{1},\dots,\psi_{j})^{\sim} \circ \bar{\sigma}_{3}.     
\end{equation}
Let us as usual consider $w = r_{i} u_{i} = v_{i} s_{i} \in R_{N+1}$. 
On the one side, $(d(\psi'_{j+1})^{\sim} \circ \bar{\sigma}_{3}) (1|w|1)$ is equal to 
\begin{align*}
   d(\psi'_{j+1})&^{\sim}(1|v_{i}|s_{i,(1)}|\bar{s}_{i,(2)}|s_{i,(3)}) 
                                                = (\psi'_{j+1})^{\sim}(v_{i}|s_{i,(1)}|\bar{s}_{i,(2)}|s_{i,(3)}) 
   \\
                                                &- (\psi'_{j+1})^{\sim}(1|v_{i}s_{i,(1)}|\bar{s}_{i,(2)}|s_{i,(3)}) 
                                                + (\psi'_{j+1})^{\sim}(1|v_{i}|s_{i}|1) 
   \\
                                                &= v_{i} \psi'_{j+1}(s_{i,(1)},\bar{s}_{i,(2)})s_{i,(3)}
                                                - \psi'_{j+1}(v_{i}s_{i,(1)},\bar{s}_{i,(2)})s_{i,(3)} 
   \\
                                                &= v_{i} \psi'_{j+1}(s_{i,(1)},\bar{s}_{i,(2)})
                                                - \psi'_{j+1}(v_{i}s_{i,(1)},\bar{s}_{i,(2)})s_{i,(3)}
   \\
                                                &= v_{i} \psi'_{j+1}(s_{i,(1)},\bar{s}_{i,(2)})
                                                - \underset{\deg(s_{i,(3)}) = 0}{\underbrace{\psi'_{j+1}(v_{i}s_{i,(1)},\bar{s}_{i,(2)})}}
                                                - \underset{\deg(s_{i,(3)}) > 0}{\underbrace{\psi'_{j+1}(v_{i}s_{i,(1)},\bar{s}_{i,(2)})s_{i,(3)}}}
   \\
                                                &= v_{i} \psi'_{j+1}(s_{i,(1)},\bar{s}_{i,(2)})
                                                - \underset{\deg(s_{i,(3)}) = 0}{\underbrace{\psi'_{j+1}(v_{i}s_{i,(1)},\bar{s}_{i,(2)})}}
                                                - \underset{\deg(s_{i,(3)}) = 1}{\underbrace{\psi'_{j+1}(v_{i}s_{i,(1)},\bar{s}_{i,(2)})s_{i,(3)}}},
\end{align*}
where we have used Fact \ref{fac:cuenta telescopica} in the second equality, that $s_{i}$ vanishes in $A$ in the third equality, and 
that $\psi'_{j+1}$ is normalized in the fourth and sixth ones. 
Moreover, since $v_{i}s_{i,(1)} \otimes \bar{s}_{i,(2)} = r_{i} \otimes u_{i}$, we see that  
\[     \psi'_{j+1}(v_{i}s_{i,(1)},\bar{s}_{i,(2)}) = \psi'_{j+1}(r_{i},u_{i}) = 0.     \]
Also, taking into account that if $\deg(s_{i,(3)}) = 1$, then 
\[     v_{i}s_{i,(1)} \otimes \bar{s}_{i,(2)} \otimes s_{i,(3)} = r_{i,(1)} \otimes \bar{r}_{i,(2)} \otimes u_{i},     \]
we obtain that $\psi'_{j+1}(v_{i}s_{i,(1)},\bar{s}_{i,(2)})s_{i,(3)} = \psi'_{j+1}(r_{i,(1)}, \bar{r}_{i,(2)}) u_{i}$, so 
\[     (d(\psi'_{j+1})^{\sim} \circ \bar{\sigma}_{3}) (1|w|1) = v_{i} \psi'_{j+1}(s_{i,(1)},\bar{s}_{i,(2)}) - \psi'_{j+1}(r_{i,(1)}, \bar{r}_{i,(2)}) u_{i}.     \]
By its very definition, $\psi'_{j}(r_{(1)},\bar{r}_{(2)}) = \varphi_{N-j}(r)$, for all $r \in R$, $j=1, \dots, N$. 
Therefore,
\begin{equation}
\label{eq:mit1}
\begin{split}
     (d(\psi'_{j+1})^{\sim} \circ \bar{\sigma}_{3}) (1|w|1) &= v_{i} \varphi_{N-j-1}(s_{i}) - \varphi_{N-j-1}(r_{i}) u_{i} 
     \\
                                               &= (1_{V} \otimes \varphi_{N-j-1} - \varphi_{N-j-1} \otimes 1_{V})(w).     
\end{split}
\end{equation}

On the other side, $(\mathrm{sq}(\psi_{1},\dots,\psi_{j})^{\sim} \circ \bar{\sigma}_{3}) (1|w|1)$ is equal to 
\begin{align*}
   \mathrm{sq}(\psi_{1}&,\dots,\psi_{j})^{\sim}(1|v_{i}|s_{i,(1)}|\bar{s}_{i,(2)}|s_{i,(3)})
   \\   
   &= \sum_{i=1}^{j} (\psi_{i}(v_{i},\psi_{j+1-i}(s_{i,(1)},\bar{s}_{i,(2)})) - \psi_{i}(\psi_{j+1-i}(v_{i},s_{i,(1)}),\bar{s}_{i,(2)}))s_{i,(3)} 
   \\
   &= \sum_{i=1}^{j} (\psi_{i}(v_{i},\psi_{j+1-i}(s_{i,(1)},\bar{s}_{i,(2)})) - \psi_{i}(\psi_{j+1-i}(v_{i},s_{i,(1)}),\bar{s}_{i,(2)}))
   \\
   &= \sum_{i=1}^{j} (\psi_{i}(v_{i},\psi_{j+1-i}(s_{i,(1)},\bar{s}_{i,(2)})) - \psi_{i}(\psi_{j+1-i}(\bar{r}_{i,(1)},r_{i,(2)}),u_{i}))
   \\
   &= \psi_{j}(v_{i},\psi_{1}(s_{i,(1)},\bar{s}_{i,(2)})) - \psi_{j}(\psi_{1}(\bar{r}_{i,(1)},r_{i,(2)}),u_{i}),
\end{align*} 
where we have used in the third equality that $\psi_{i}$ is normalized. 
In the last one we have used the following simple fact which we have already explained: 
$\psi_{j+1-i}$ has degree $-(j+1-i) = i-j-1$, so 
$\psi_{j+1-i}$ applied to an element of degree $N$ gives an element of degree $N+i-j-1$. 
In consequence, the elements  $v_{i} \otimes \psi_{j+1-i}(s_{i,(1)},\bar{s}_{i,(2)})$ and $\psi_{j+1-i}(\bar{r}_{i,(1)},r_{i,(2)}) \otimes u_{i}$ have degree $N+j-i \leq N$ 
(for $i \leq j$) and the degree is exactly $N$ if and only if $i=j$. 
The normalization of $\psi_{j}$ forces only to consider the terms with degree $N$ and the last equality follows. 

By its very definition, $\psi_{1}(r_{(1)},\bar{r}_{(2)}) = \psi'_{1}(r_{(1)},\bar{r}_{(2)}) = \varphi_{N-1}(r)$, for $r \in R$, so  
\[     (\mathrm{sq}(\psi_{1},\dots,\psi_{j})^{\sim} \circ \bar{\sigma}_{3}) (1|w|1) = \psi_{j}(v_{i},\varphi_{N-1}(s_{i})) - \psi_{j}(\varphi_{N-1}(r_{i}),u_{i}).     \]
Thus, $(\mathrm{sq}(\psi_{1},\dots,\psi_{j})^{\sim} \circ \bar{\sigma}_{3}) (1|w|1)$ is given by 
\begin{align*}
   \psi_{j}(v_{i}\otimes \varphi_{N-1}(s_{i}) - &\varphi_{N-1}(r_{i}) \otimes u_{i})
   \\
   &= (\varphi_{N-j}^{\sim} \circ \bar{\tau}_{2} + \eta_{j}^{\sim})(1|\underset{\phantom{r} \in R}{\underbrace{v_{i}\otimes \varphi_{N-1}(s_{i}) - \varphi_{N-1}(r_{i}) \otimes u_{i}}}|1)
   \\
   &= (\varphi_{N-j}^{\sim} \circ \bar{\tau}_{2})(1|\underset{\phantom{r} \in A^{\otimes 2}}{\underbrace{v_{i}\otimes \varphi_{N-1}(s_{i}) - \varphi_{N-1}(r_{i}) \otimes u_{i}}}|1)
   \\
   &= (\varphi_{N-j}^{\sim} (1|\underset{\phantom{r} \in R}{\underbrace{v_{i}\otimes \varphi_{N-1}(s_{i}) - \varphi_{N-1}(r_{i}) \otimes u_{i}}}|1)
   \\
   &= \varphi_{N-j} \circ (1_{V} \otimes \varphi_{N-1} - \varphi_{N-1} \otimes 1_{V})(w),
\end{align*} 
where we have used that, by identity \eqref{eq:pbwfi1}, $v_{i}\otimes \varphi_{N-1}(s_{i}) - \varphi_{N-1}(r_{i}) \otimes u_{i}$ can be seen as 
a relation decomposition of an element $r$ of $R$ in the second equality, and that $\eta_{j}^{\sim}$ vanishes over it in the third equality. 
Finally, in the penultimate equality we have used that 
\[     \bar{\tau}_{2}(1|v_{i}\otimes \varphi_{N-1}(s_{i}) - \varphi_{N-1}(r_{i}) \otimes u_{i}|1) = 1|(v_{i}\otimes \varphi_{N-1}(s_{i}) - \varphi_{N-1}(r_{i}) \otimes u_{i})|1,     \]
where in the first member $v_{i}\otimes \varphi_{N-1}(s_{i}) - \varphi_{N-1}(r_{i}) \otimes u_{i}$ is seen as an element of $A^{\otimes 2}$, 
whereas in the second one it is regarded as an element of $R$. 
This thus implies that 
\begin{equation}
\label{eq:mit2}
     (\mathrm{sq}(\psi_{1},\dots,\psi_{j})^{\sim} \circ \bar{\sigma}_{3}) (1|w|1) = \varphi_{N-j} \circ (1_{V} \circ \varphi_{N-1} - \varphi_{N-1} \otimes 1_{V})(w).   
\end{equation}

Now, since \eqref{eq:pbwfi1}, \eqref{eq:pbwfi2} and \eqref{eq:pbwfi3} hold, and putting together \eqref{eq:mit1} and \eqref{eq:mit2}, 
we see that equation \eqref{eq:Gdefs3} holds. 
The standard identity \eqref{eq:sn} for the homotopy $s_{\bullet}$ tells us that 
\begin{align*}
       d(\psi'_{j+1})^{\sim} \circ (1_{\bar{C}_{3}(A)} - \bar{\sigma}_{3} \circ \bar{\tau}_{3}) 
       &= d(\psi'_{j+1})^{\sim} \circ (\bar{b}_{4} \circ s_{3} + s_{2} \circ \bar{b}_{3})
       \\
       &= d(\psi'_{j+1})^{\sim} \circ s_{2} \circ \bar{b}_{3} 
       = d(d(\psi'_{j+1})^{\sim} \circ s_{2}),     
\end{align*}
so 
\[     d(\psi'_{j+1})^{\sim} \circ \bar{\sigma}_{3} \circ \bar{\tau}_{3} = d(\psi'_{j+1})^{\sim} - d(d(\psi'_{j+1})^{\sim} \circ s_{2}).     \]
Also, we see that 
\begin{align*}
     \mathrm{sq}(\psi_{1},\dots,\psi_{j})^{\sim} \circ (1_{\bar{C}_{3}(A)} - \bar{\sigma}_{3} \circ \bar{\tau}_{3}) 
       &= \mathrm{sq}(\psi_{1},\dots,\psi_{j})^{\sim} \circ (\bar{b}_{4} \circ s_{3} + s_{2} \circ \bar{b}_{3}) 
       \\
       &= \mathrm{sq}(\psi_{1},\dots,\psi_{j})^{\sim} \circ s_{2} \circ \bar{b}_{3}
       \\
       &= d(\mathrm{sq}(\psi_{1},\dots,\psi_{j})^{\sim} \circ s_{2}),     
\end{align*}       
\textit{i.e.} 
\[     \mathrm{sq}(\psi_{1},\dots,\psi_{j})^{\sim} \circ \bar{\sigma}_{3} \circ \bar{\tau}_{3} = \mathrm{sq}(\psi_{1},\dots,\psi_{j})^{\sim} - d(\mathrm{sq}(\psi_{1},\dots,\psi_{j})^{\sim} \circ s_{2}).     \]
Now, equation \eqref{eq:Gdefs3} yields that 
\[     -d((\psi'_{j+1})^{\sim} - (d(\psi'_{j+1})^{\sim} \circ s_{2}+ \mathrm{sq}(\psi_{1},\dots,\psi_{j})^{\sim} \circ s_{2})) 
       = \mathrm{sq}(\psi_{1},\dots,\psi_{j})^{\sim}.     \]
We define $\eta_{j+1}^{\sim} = d(\psi'_{j+1})^{\sim} \circ s_{2} + \mathrm{sq}(\psi_{1},\dots,\psi_{j})^{\sim} \circ s_{2}$ 
and $\psi_{j+1}^{\sim} = (\psi'_{j+1})^{\sim} - \eta_{j+1}^{\sim}$. 
It is easy to see that $\eta_{j+1}^{\sim}$ has degree $-j-1$. 
We just need to prove that it is extranormalized in order to end this recursive process, since $\psi_{j+1}^{\sim}$ defines a continuation of the 
$j$-th level deformation defined by $\psi_{1}^{\sim}, \dots, \psi_{j}^{\sim}$. 

Let us prove that $\eta_{j+1}^{\sim}$ is extranormalized. 
For this, consider $1|a_{i}|b_{i}|1 \in \bar{C}_{2}(A)$ to be normalized or a relation decomposition. 
Then, $(d(\psi'_{j+1})^{\sim} \circ s_{2}) (1|a_{i}|b_{i}|1)$ is equal to 
\begin{align*}
       d(\psi'_{j+1})^{\sim} &(1|a_{i}|b_{i,(1)}|\bar{b}_{i,(2)}|b_{i,(3)}) 
       \\
       &= (\psi'_{j+1})^{\sim} (a_{i}|b_{i,(1)}|\bar{b}_{i,(2)}|b_{i,(3)} - 1|a_{i}b_{i,(1)}|\bar{b}_{i,(2)}|b_{i,(3)} 
       + 1|a_{i}|b_{i}|1) 
       \\
       &= a_{i}\psi'_{j+1}(b_{i,(1)},\bar{b}_{i,(2)})b_{i,(3)} - \psi'_{j+1}(a_{i}b_{i,(1)},\bar{b}_{i,(2)})b_{i,(3)} 
       + \psi'_{j+1}(a_{i},b_{i}) 
       \\
       &=  - \psi'_{j+1}(a_{i}b_{i,(1)},\bar{b}_{i,(2)}) + \psi'_{j+1}(a_{i},b_{i}) 
       \\
       &=  - \varphi_{N-j-1} (\bar{\tau}_{2}(1|a_{i}b_{i,(1)}|\bar{b}_{i,(2)}|1) - \bar{\tau}_{2}(1|a_{i}|b_{i}|1)), 
\end{align*}
where we have used Fact \ref{fac:cuenta telescopica} in the first equality and the normalization of $\psi'_{j+1}$ on the fourth equality. 
By definition, we see that $\bar{\tau}_{2}(1|a_{i}b_{i,(1)}|\bar{b}_{i,(2)}|1) = \bar{\tau}_{2}(1|a_{i}|b_{i}|1) = 0$ if $1|a_{i}|b_{i}|1$ is normalized and 
that $\bar{\tau}_{2}(1|a_{i}b_{i,(1)}|\bar{b}_{i,(2)}|1) = \bar{\tau}_{2}(1|a_{i}|b_{i}|1)$ if $1|a_{i}|b_{i}|1$ is a relation decomposition. 
In any case, we conclude that $d(\psi'_{j+1})^{\sim} \circ s_{2} (1|a_{i}|b_{i}|1)$ vanishes. 
       
On the other hand, we get
\begin{multline*}
     (\mathrm{sq}(\psi_{1},\dots,\psi_{j})^{\sim} \circ s_{2}) (1|a_{i}|b_{i}|1) 
     = \mathrm{sq}(\psi_{1},\dots,\psi_{j})^{\sim} (1|a_{i}|b_{i,(1)}|\bar{b}_{i,(2)}|b_{i,(3)}) 
     \\
     = \sum_{i=1}^{j}\Big((\psi_{i}(a_{i},\psi_{j+1-i}(b_{i,(1)},\bar{b}_{i,(2)})) - \psi_{i}(\psi_{j+1-i}(a_{i},b_{i,(1)}),\bar{b}_{i,(2)}))b_{i,(3)}\Big),  
\end{multline*}
which vanishes by the normalization of $\psi_{j+1-i}$.
Hence, $\eta_{j+1}^{\sim}$ is extranormalized. 

We may summarize the previous results as follows.
\begin{proposition}
\label{prop:pbwtoger}
Let $U = TV/\cl{P}$ be a filtered algebra, with $P \subset F^{N}$, and let $A = TV/\cl{R}$, with $R = \pi_{N}(P) \subset V^{\otimes N}$, be the corresponding $N$-homogeneous algebra, 
which we assume to satisfy that $\Tor_{A}^{3}(k,k)$ is concentrated in degree $N+1$. 
We assume that $U$ is a weak PBW-deformation of $A$, 
\textit{i.e.} that \eqref{eq:pbw1} holds and that the $k^{e}$-linear maps $\varphi_{j} : R \rightarrow V^{\otimes j}$ (for $j=0,\dots,N-1$) 
which describe $P$ out from $R$ satisfy \eqref{eq:pbwfi1}, \eqref{eq:pbwfi2} and \eqref{eq:pbwfi3}.  
Set $(\psi'_{j})^{\sim} = \varphi_{N-j} \circ \bar{\tau}_{2}$ for $j=1, \dots, N$ and zero otherwise. 
We define $\psi_{j}^{\sim} : \bar{C}_{2}(A) \rightarrow A$ for $j \in \NN$ recursively. 
First, $\psi_{1}^{\sim} = (\psi'_{1})^{\sim}$. 
For $j \in \NN$, after having defined $\psi_{1},\dots,\psi_{j}$, we set $\eta_{j+1}^{\sim} = d(\psi'_{j+1})^{\sim} \circ s_{2} + \mathrm{sq}(\psi_{1},\dots,\psi_{j})^{\sim} \circ s_{2}$ 
and $\psi_{j+1}^{\sim} = (\psi'_{j+1})^{\sim} - \eta_{j+1}^{\sim}$. 
Then the $\psi_{\bullet}$ are normalized morphisms that define a graded deformation of $A$. 
\end{proposition}

\section{Main theorems}
\label{sec:main}

The following theorems provide a description of the previous constructions at the level of algebras, which contains the one given in Sec. 4.6 in \cite{BGa96}, 
where the authors explored only one direction under the assumption of $k = F$ a field. 
Moreover, we also prove that conditions \eqref{eq:pbw1} and \eqref{eq:pbw2} 
are equivalent to the fact that $U$ satisfies the PBW property (\textit{cf.}~Thm. 4.1 of \cite{BGa96} and Thm. 3.4 of \cite{BG06}) . 

\begin{theorem}
\label{teo:pbw}
Let $A$ be an $N$-homogeneous algebra satisfying that $\Tor_{3}^{A}(k,k)$ is concentrated in degree $N+1$. 
Let us consider a graded deformation $A_{t}$ of a $A$, which we suppose to be given by 
normalized maps $\{ \psi_{j} : (A/k)^{\otimes 2} \rightarrow A \}_{j \in \NN}$.  
We apply Proposition \ref{prop:gertopbw} to produce maps $\{ \varphi_{j} : R \rightarrow V^{\otimes j} \}_{0 \leq j < N }$, and 
to obtain thus a filtered algebra $U = TV/\cl{P}$, with $P = \{ r - \sum_{j=0}^{N-1} \varphi_{j}(r) : r \in R \}$. 
Then, there exists an isomorphism of filtered $k$-algebras $U \rightarrow A_{t}/\cl{t-1}$.

Conversely, let us consider a filtered algebra $U = TV/\cl{P}$, with $P \subset F^{N}$ and let $A = TV/\cl{R}$, with $R = \pi_{N}(P) \subset V^{\otimes N}$, 
be the corresponding $N$-homogeneous algebra, which we assume to satisfy that $\Tor_{A}^{3}(k,k)$ is concentrated in degree $N+1$. 
We assume that $U$ is a weak PBW-deformation of $A$, and define a deformation $A_{t}$ of $A$ following Proposition \ref{prop:pbwtoger}. 
We again see that there exists an isomorphism of filtered $k$-algebras $U \rightarrow A_{t}/\cl{t-1}$. 

In both cases, $U$ is a PBW-deformation of $A$ and the induced morphism 
\[     A \overset{p}{\rightarrow} \mathrm{gr}(U) \rightarrow \mathrm{gr}(A_{t}/\cl{t-1}) \overset{\rho}{\rightarrow} A     \] 
is the identity. 
\end{theorem}
\noindent\textbf{Proof.}
Since the proof is similar for both implications, except for some minor changes, we shall only treat each case separately when necessary.

We consider the $k^{e}$-linear map $\mathrm{inc} : V \rightarrow A_{t}/\cl{t-1}$, given by the composition of 
the inclusion $V \rightarrow A$, the canonical map $A \rightarrow A_{t}$ and the projection $A_{t} \rightarrow A_{t}/\cl{t-1}$. 
This induces a morphism of $k$-algebras $q' : TV \rightarrow A_{t}/\cl{t-1}$. 
It is clear that $q'$ respects the filtrations (where we recall that $TV$ is filtered by $F^{\bullet}$). 

We shall see that $q'(P) = 0$. 
In order to prove this statement, take $r \in R$ and let $r - \sum_{i=0}^{N-1} \varphi_{i}(r) \in P$ be a generic element. 
We will show that $q'$ vanishes over it, using that the associated maps $\psi_{\bullet}$ are normalized.  

First, we remark the fact that, if $\alpha \in V^{\otimes j}$, with $j < N$, then $q'(\alpha) = \alpha$, which can be proved as follows. 
It suffices to treat the case $\alpha = v_{1} \dots v_{j}$, for $v_{1}, \dots, v_{j} \in V$. 
By definition $q'$ acts as the identity for $j = 0, 1$. 
Let us thus assume that $j \geq 2$ and prove the statement by induction on $j$. 
We assume that $q'$ acts as the identity on $F^{j-1}$ and we shall prove that it does the same on $F^{j}$ ($j < N$). 
The inductive hypothesis implies that $q'(v_{1} \dots v_{j-1}) = v_{1} \times \dots \times v_{j-1} = v_{1} \dots v_{j-1}$. 
Hence,
\begin{align*}
    q'(v_{1} \dots v_{j}) &= v_{1} \times \dots \times v_{j} = (v_{1} \times \dots \times v_{j-1}) \times v_{j} = (v_{1} \dots v_{j-1}) \times v_{j}
    \\
    & = v_{1} \dots v_{j-1} v_{j} + \sum_{l \geq 1} \psi_{l}(v_{1} \dots v_{j-1},v_{j})                                  
    = v_{1} \dots v_{j-1} v_{j},     
\end{align*}
since the maps $\psi_{\bullet}$ are normalized. 

Now, take $r = v_{1,i} \dots v_{N,i} \in R$ (summation understood). 
Then, 
\begin{align*}
     q'(r) &= v_{1,i} \times \dots \times v_{N,i} = (v_{1,i} \times \dots \times v_{N-1,i}) \times v_{N,i} = (v_{1,i} \dots v_{N-1,i}) \times v_{N,i} 
     \\        
          &= v_{1,i} \dots v_{N-1,i} v_{N,i} + \sum_{l \in \NN} \psi_{l}(v_{1,i} \dots v_{N-1,i} ,v_{N,i}).     
\end{align*}
If we are considering the first statement, using that, by definition, $\varphi_{N-l}(r) = \psi_{l} (r_{1},\bar{r}_{2})$, for $l= 1, \dots, N$, 
and by degree reasons $\psi_{l} (r_{1},\bar{r}_{2})=0$ for $l > N$, 
then $q'(r) = r + \sum_{l=1}^{N} \varphi_{N-l}(r) = \sum_{l=0}^{N-1} \varphi_{l}(r)$, since $r$ vanishes in $A$ and $t=1$ in $A_{t}/\cl{t-1}$. 
For the second statement, since $v_{1,i} \dots v_{N-1,i} \otimes v_{N,i} \in A^{\otimes 2}$ is a relation decomposition 
we see that $\eta_{j}(v_{1,i} \dots v_{N-1,i} , v_{N,i})$ vanishes, 
because it is extranormalized, and $\bar{\tau}_{2}(1|v_{1,i} \dots v_{N-1,i} | v_{N,i}|1) = 1|r|1$. 
Hence, 
\[     \psi_{l}(v_{1,i} \dots v_{N-1,i} ,v_{N,i}) = \psi'_{l}(v_{1,i} \dots v_{N-1,i} ,v_{N,i}) = \begin{cases} 
                                                                                                      \varphi_{N-l}(r), &\text{if $l = 1, \dots, N$,}
                                                                                                      \\     
                                                                                                      0,  &\text{else.}
                                                                                                                          \end{cases}
\]
As a consequence, we again have that $q'(r) = r + \sum_{l=1}^{N} \varphi_{N-l}(r) = \sum_{l=0}^{N-1} \varphi_{l}(r)$, since $r$ vanishes in $A$ and $t=1$ in $A_{t}/\cl{t-1}$. 
Therefore, $q'(r - \sum_{i=0}^{N-1} \varphi_{i}(r)) = 0$ and $q'$ thus induces a morphism of filtered algebras $q : U \rightarrow A_{t}/\cl{t-1}$. 

We will now prove that $q$ is an isomorphism. 
Since 
\[     A \overset{p}{\rightarrow} \mathrm{gr}(U) \overset{\mathrm{gr}(q)}{\rightarrow} \mathrm{gr}(A_{t}/\cl{t-1}) \overset{\rho}{\rightarrow} A     \]
is the identity, $\rho$ is an isomorphism and $p$ is surjective, we conclude that $p$ is an isomorphism and $\mathrm{gr}(q)$ is also an isomorphism. 
Hence, $q$ is an isomorphism and $U$ is a PBW-deformation of A.
\qed

\begin{remark}
The preceding theorem implies that $A_{t}/\cl{t-1}$ may be seen as a PBW-deformation of $A$ equivalent to $U$, where the morphism 
$A \rightarrow \mathrm{gr}(A_{t}/\cl{t-1})$ is the inverse of $\rho$. 
Furthermore, by Lemma \ref{lema:sim} we see that the procedure in Subsection \ref{subsec:gertopbw} may be defined in the set of equivalence classes of deformations 
and it sends equivalent deformations to equivalent PBW-deformations. 
So it defines a map
\[     \mathrm{gp} : \{ \text{eq. classes of deformations of $A$} \} \rightarrow \{ \text{eq. classes of PBW-deformations of $A$} \}.     \]
\end{remark}

\begin{theorem}
\label{teo:pbw+}
Let $A$ be an $N$-homogeneous algebra satisfying that $\Tor_{3}^{A}(k,k)$ is concentrated in degree $N+1$.
Let us consider a graded deformation $A_{t}$ of a $A$ by normalized maps $\{ \psi_{j} : (A/k)^{\otimes 2} \rightarrow A \}_{j \in \NN}$.  
We apply the construction of Proposition \ref{prop:gertopbw} to produce maps $\{ \varphi_{j} : R \rightarrow V^{\otimes j} \}_{0 \leq j < N }$, and 
to obtain thus a filtered algebra $U = TV/\cl{P}$, with $P = \{ r - \sum_{j=0}^{N} \varphi_{j}(r) : r \in R \}$. 
Then, there exists an isomorphism of graded $k[t]$-algebras $R(U) \rightarrow A_{t}$ such that the induced morphism
$A \simeq R(U)/\cl{t} \rightarrow A_{t}/\cl{t} \simeq A$ is the identity, where $A \simeq R(U)/\cl{t}$ is the map described in Subsection \ref{subsec:ger}. 

Conversely, let us consider a filtered algebra $U = TV/\cl{P}$, for $P \subset F^{N}$ and let $A = TV/\cl{R}$, for $R = \pi_{N}(P) \subset V^{\otimes N}$, 
be the corresponding $N$-homogeneous algebra, which we assume to satisfy that $\Tor_{A}^{3}(k,k)$ is concentrated in degree $N+1$. 
We suppose that $U$ is a weak PBW-deformation of $A$, and define a deformation $A_{t}$ of $A$ following Proposition \ref{prop:pbwtoger}. 
We again see that there exists an isomorphism of graded $k[t]$-algebras $R(U) \rightarrow A_{t}$ such that the induced morphism
$A \simeq R(U)/\cl{t} \rightarrow A_{t}/\cl{t} \simeq A$ is the identity, where $A \simeq R(U)/\cl{t}$ is the map described in Subsection \ref{subsec:ger}. 
\end{theorem}
\noindent\textbf{Proof.} 
From Proposition \ref{prop:fibra}, we see that, given any graded deformation $A_{t}$ of $A$, there exists an isomorphism of graded $k[t]$-algebras 
$\mathrm{com} : R(A_{t}/\cl{t-1}) \rightarrow A_{t}$. 
Since $q : U \rightarrow A_{t}/\cl{t-1}$ is an isomorphism of filtered algebras, $R(q) : R(U) \rightarrow R(A_{t}/\cl{t-1})$ is an isomorphism of graded $k[t]$-algebras.  
The composition $\mathrm{com} \circ R(q)$ gives the desired isomorphism. 
Using the canonical isomorphism $R(\place)/\cl{t} \simeq \mathrm{gr}(\place)$, we see that 
$A \simeq R(U)/\cl{t} \rightarrow R(A_{t}/\cl{t-1}) \simeq A$ identifies with 
\[     A \overset{p}{\rightarrow} \mathrm{gr}(U) \overset{\mathrm{gr}(q)}{\rightarrow} \mathrm{gr}(A_{t}/\cl{t-1}) \overset{\rho}{\rightarrow} A,     \]
which we have already seen to be the identity. 
\qed

\begin{remark}
The preceding theorem implies that $R(U)$ may be seen as a deformation of $A$ equivalent to $A_{t}$. 
Moreover, it also says that the procedure in Subsection \ref{subsec:pbwtoger} sends equivalence classes of PBW-deformations to equivalent deformations, \textit{i.e.} 
it defines a map
\[     \mathrm{pg} : \{ \text{eq. classes of PBW-deformations of $A$} \} \rightarrow \{ \text{eq. classes of deformations of $A$} \}.     \]
\end{remark}

\begin{theorem}
\label{teo:final}
By the previous theorems we see that the procedure performed in Subsection \ref{subsec:gertopbw} consists 
in $A_{t} \mapsto A_{t}/\cl{t-1}$ whereas the one done in Subsection \ref{subsec:pbwtoger} is just $U \rightarrow R(U)$, both up to equivalence. 
Thus Proposition \ref{prop:fibra}, Theorem \ref{teo:pbw} and Theorem \ref{teo:pbw+} imply that $\mathrm{pg}$ and $\mathrm{gp}$ are mutually inverse. 
\end{theorem}

We finish by exhibiting two examples. 
\begin{example}
Let $k = F$ be a field, $V = \mathrm{span}_{k}\cl{x,y}$ a $k$-vector space of dimension $2$, and $A = TV/\cl{[x,y]}$ a quadratic algebra, 
so $R \subset V^{\otimes 2}$ has a basis formed by a unique element $r = [x,y]$. 
We note that $R_{3} = 0$. 

We can provide in this case the maps $\bar{\tau}_{\bullet}$, for $\bullet = 1, 2$, in their complete domain of definition:
\begin{align*}
 \bar{\tau}_{1} (1|y^{n} x^{m}|1) 
   &= \sum_{\alpha_{1} + \alpha_{2} = n - 1} y^{\alpha_{1}} |y|y^{\alpha_{2}} x^{m} + \sum_{\beta_{1} + \beta_{2} = m - 1} 
                    y^{n} x^{\beta_{1}}|x|x^{\beta_{2}},
 \\
 \bar{\tau}_{2} (1|y^{n_{1}} x^{m_{1}}|y^{n_{2}} x^{m_{2}}|1) 
   &= \sum_{\alpha_{1} + \alpha_{2} = n_{2} - 1} \sum_{\beta_{1} + \beta_{2} = m_{1} - 1} y^{n_{1}+\alpha_{1}} x^{\beta_{1}}|r|y^{\alpha_{2}} x^{m_{2}+\beta_{2}}.
\end{align*}
Furthermore, the homotopies $s_{\bullet}$, for $\bullet = 1, 2$, are given as follows. 
The image of $1|y^{n} x^{m}|1$ under $s_{1}$ is the class in $\bar{C}_{2}(A)$ of 
\[      
        - \sum_{\alpha + \alpha_{2} = n - 1} 1|y^{\alpha_{1}} |y|y^{\alpha_{2}} x^{m} - \sum_{\beta_{1} + \beta_{2} = m - 1} 
                    1|y^{n} x^{\beta_{1}}|x|x^{\beta_{2}},     \]
and the image of $1|y^{n_{1}} x^{m_{1}}|y^{n_{2}} x^{m_{2}}|1$ under $s_{2}$ is the class in $\bar{C}_{3}(A)$ of 
\begin{multline*}
\sum_{\begin{scriptsize}\begin{matrix}\alpha_{1} + \alpha_{2} = n_{2} - 1\\\beta_{1} + \beta_{2} = m_{1} - 1\end{matrix}\end{scriptsize}} 
   \Big(1|y^{n_{1}+\alpha_{1}} x^{\beta_{1}}|y|x|y^{\alpha_{2}} x^{m_{2}+\beta_{2}} 
  - 1|y^{n_{1}+\alpha_{1}} x^{\beta_{1}}|x|y|y^{\alpha_{2}} x^{m_{2}+\beta_{2}} \Big)
 \\ 
   + \sum_{\alpha_{1} + \alpha_{2} = n_{2} - 1} 1|y^{n_{1}} x^{m_{1}}|y^{\alpha_{1}}|y|y^{\alpha_{2}} x^{m_{2}} 
   + \sum_{\beta_{1} + \beta_{2} = m_{2} - 1} 1|y^{n_{1}} x^{m_{1}}| y^{n_{2}} x^{\beta_{1}}|x|x^{\beta_{2}}.
\end{multline*}

Suppose that $U = TV/\cl{[x,y]-y}$. 
Note that the maps defining the filtration are $\varphi_{0} = 0$ and $\varphi_{1}(r) = y$, which obviously satisfy the weak PBW-property. 
A simple inductive argument then implies that the cochains defining the deformed product are 
\[     \psi_{l}(y^{n_{1}} x^{m_{1}},y^{n_{2}} x^{m_{2}}) = \left( \begin{matrix} m_{1} \\ l \end{matrix} \right) n_{2}^{l} y^{n_{1}+n_{2}} x^{m_{1}+m_{2}-l}.      \]

On the other hand, if $U = TV/\cl{[x,y]-1}$, the so-called \emph{Weyl algebra}, the maps defining the filtration are $\varphi_{0}(r) = 1$ and $\varphi_{1} = 0$, which also satisfy the weak PBW-property. 
It is direct to prove in this case that the cochains giving the deformed product are of the form 
\[     \psi_{l} (y^{n_{1}} x^{m_{1}},y^{n_{2}} x^{m_{2}}) = \begin{cases}   
                                                              i! \left( \begin{matrix} m_{1} \\ i \end{matrix} \right)  
                                                                   \left( \begin{matrix} n_{2} \\ i \end{matrix} \right) y^{n_{1}+n_{2}-i} x^{m_{1}+m_{2}-i}, &\text{if $l = 2 i$, $i \in \NN$,}
                                                              \\
                                                              0, &\text{else.}      
                                                            \end{cases}
\]
\end{example}

\begin{example}
Let $k = F$ be a field, $V = \mathrm{span}_{k}\cl{x}$ a $k$-vector space of dimension $1$, and $A = TV/\cl{x^{N}}$ an $N$-homogeneous algebra, 
so $R \subset V^{\otimes N}$ has a basis given by $r = x^{N}$. 
Notice that $R_{3} = \mathrm{span}_{k}\cl{x^{N+1}}$. 

In this case, the comparison map $\bar{\tau}_{1}$ is actually given in \eqref{eq:tau1} in its complete domain of definition $\bar{C}_{1}(A)$. 
As for $\bar{\tau}_{2}$, its full expression may be given by (see also \cite{BACH}):
\[      \bar{\tau}_{2} (1|x^{m_{1}}|x^{m_{2}}|1) = \begin{cases}
           1|x^{N}|x^{m_{2}+m_{1}-N}, &\text{if $m_{1} + m_{2} \geq N$},
           \\
           0, &\text{else},
          \end{cases}
\]
where we consider $1 \leq m_{1}, m_{2} \leq N-1$. 
Furthermore, the homotopy $s_{1}$ is also given in \eqref{eq:s1} in its complete domain of definition. 
On the other hand, the homotopy $s_{2}$ is given by the obvious extension  
\[        s_{2} (1|x^{m_{1}}|x^{m_{2}}|1) = \sum_{i=1}^{m_{2}-i-1} 1|x^{m_{1}}|x^{i}|x|x^{m_{2}-i-1}, 
\]
where $1 \leq m_{1}, m_{2} \leq N-1$, and, as previously stated, we suppose that the right member lies in $\bar{C}_{3}(A)$. 

Let $f = \sum_{i=0}^{N-1} a_{i} x^{i} \in k[x]$ be a polynomial of degree less than or equal to $N-1$, and suppose that $U = TV/\cl{x^{N}-f}$. 
As usual we assume that $a_{i} = 0$, if $i \notin \{ 0, \dots, N-1\}$. 
Note that the maps defining the filtration are $\varphi_{i}(r) = a_{i} x^{i}$, for $i = 0, \dots, N-1$, which obviously satisfy the weak PBW-property. 
A simple inductive argument then implies that the cochains defining the deformed product are 
\[     \psi_{l}(x^{m_{1}},x^{m_{2}}) = \begin{cases}   
                                                              a_{N-l} x^{m_{1}+m_{2}-l}, &\text{if $m_{1} + m_{2} \geq N$,}
                                                              \\
                                                              0, &\text{if $m_{1}+m_{2} < N$,}      
                                                            \end{cases},      \]
as one could have expected. 
\end{example}



\medskip

\noindent Estanislao Herscovich (corresponding author),
\\
Fakult\"at f\"ur Mathematik,
\\
Universit\"at Bielefeld,
\\
D-33615 Bielefeld,
\\
Germany,
\\
\href{mailto:eherscov@math.uni-bielefeld.de}{eherscov@math.uni-bielefeld.de}

\medskip

\noindent Andrea Solotar,
\\
Departamento de Matem\'atica, FCEyN,
\\
Instituto de Matem\'atica Luis Santal\'o, IMAS-CONICET,
\\
Universidad de Buenos Aires,
\\
1428, Buenos Aires,
\\
Argentina,
\\
\href{mailto:asolotar@dm.una.ar}{asolotar@dm.una.ar}

\medskip

\noindent Mariano Su\'arez-\'Alvarez,
\\
Departamento de Matem\'atica, FCEyN,
\\
Universidad de Buenos Aires,
\\
1428, Buenos Aires,
\\
Argentina,
\\
\href{mailto:mariano@dm.uba.ar}{mariano@dm.uba.ar}

\end{document}